\documentclass[amsppt,11pt]{amsart}
\usepackage{curves}
\usepackage{amsmath}

\newtheorem {theorem}{Theorem}[section]

\newtheorem {proposition}[theorem]{Proposition}
\newtheorem {lemma}[theorem]{Lemma}
\newtheorem {corollary}[theorem]{Corollary}

\newcounter{conjecture}
\newcounter{remark}

\newenvironment{remark}{\medskip{\bf Remark
\theremark.}
\addtocounter{remark}{1}}{}
\newcommand{\eqnsection}{
          \renewcommand{\theequation}{\thesection.\arabic{equation}}
           \makeatletter
           \csname @addtoreset\endcsname{equation}{section}
           \makeatother}

\newcommand{\be}{{\begin{equation}}}
\newcommand{\ee}{{\end{equation}}}
\def \bt{\begin{theorem}}
\def \et{\end{theorem}}
\def \bea{\begin{eqnarray}}
\def \eea{\end{eqnarray}}
\def \bas{\begin{eqnarray*}}
\def \eas{\end{eqnarray*}}

\def \al{\alpha}
\def \bb{\beta}
\def \ga{\gamma}
\def \Ga{\Gamma}
\def \de{\delta}

\def \ep{\epsilon}

\newcommand{\eps}{\varepsilon}

\def \la{\lambda}

\def \om{\omega}
\def \Om{\Omega}

\def \si{\sigma}

\def \th{\theta}

\def \ze{\zeta}

\def \ff{\infty}

\def \wt{\widetilde}

\def \rar{\rightarrow}

\newcommand{\ls}[1]
           {\dimen0=\fontdimen6\the\font \lineskip=#1\dimen0
\advance\lineskip.5\fontdimen5\the\font
\advance\lineskip-\dimen0
\lineskiplimit=.9\lineskip \baselineskip=\lineskip
\advance\baselineskip\dimen0 \normallineskip\lineskip
\normallineskiplimit\lineskiplimit
\normalbaselineskip\baselineskip
\ignorespaces }
\newcommand{\req}[1]{(\ref{#1})}

\def \Z{{\Bbb{Z}}}

\def \FF{{\mathcal F}}
\def \GG{{\mathcal G}}
\def \HH{{\mathcal H}}

\def \({\left(}
\def \){\right)}

\def \lc{\left\{}
\def \rc{\right\}}

\def \nn{\nonumber}

\def \bc{\begin{center} }
\def \ec{\end{center} }
\begin{document}

\eqnsection
\newcommand{\Ini}{{I_{n,i}}}
\newcommand{\reals}{{\Bbb{R}}}
\newcommand{\F}{{\mathcal F}}
\newcommand{\D}{{\mathcal D}}
\newcommand{\Fn}{{{\mathcal F}_n}}
\newcommand{\Gn}{{{\mathcal G}_n}}
\newcommand{\Hn}{{{\mathcal H}_n}}
\newcommand{\Fp}{{{\mathcal F}^p}}
\newcommand{\Gp}{{{\mathcal G}^p}}
\newcommand{\PPP}{{\Bbb P}}
\newcommand{\Pop}{{P\otimes \PPP}}
\newcommand{\hm}{\HH^\varphi}
\newcommand{\nuw}{{\nu^W}}
\newcommand{\ths}{{\theta^*}}
\newcommand{\beq}[1]{\begin{equation}\label{#1}}
\newcommand{\eeq}{\end{equation}}
\newcommand{\integers}{{\rm I\!N}}
\newcommand{\E}{{\Bbb E}}
\newcommand{\te}{{\tilde{\delta}}}
\newcommand{\tI}{{\tilde{I}}}
\newcommand{\loge}{{\log(1/\ep)}}
\newcommand{\logen}{{\log(1/\ep_n)}}
\newcommand{\epn}{{\ep_n}}
\def\var{{\rm Var}}
\def\cov{{\rm Cov}}
\def\one{{\bf 1}}
\def\leb{{\mathcal L}eb}
\def\Ho{{\mbox{\sf H\"older}}}  
\def\thi{{\mbox{\sf Thick}}}
\def\cthi{{\mbox{\sf CThick}}}
\def\late{{\mbox{\sf Late}}}
\def\clate{{\mbox{\sf CLate}}}
\newcommand{\ffrac}[2]
          {\left( \frac{#1}{#2} \right)}
\newcommand{\calF}{{\mathcal F}}
\newcommand{\dfn}{\stackrel{\triangle}{=}}
\newcommand{\beqn}[1]{\begin{eqnarray}\label{#1}}
\newcommand{\eeqn}{\end{eqnarray}}
\newcommand{\oo}{\overline}
\newcommand{\uu}{\underline}
\newcommand{\bfcdot}{{\mbox{\boldmath$\cdot$}}}
\newcommand{\Var}{{\rm \,Var\,}}
\def\squarebox#1{\hbox to #1{\hfill\vbox to #1{\vfill}}}
\renewcommand{\qed}{\hspace*{\fill}

\vbox{\hrule\hbox{\vrule\squarebox{.667em}\vrule}\hrule}\smallskip}
\newcommand{\half}{\frac{1}{2}\:}
\newcommand{\beaa}{\begin{eqnarray*}}
\newcommand{\eeaa}{\end{eqnarray*}}
\newcommand{\calK}{{\mathcal K}}
\def\dimm{{\overline{{\rm dim}}_{_{\rm M}}}}
\def\dimp{\dim_{_{\rm P}}}
\def\htaum{{\hat\tau}_m}
\def\htaumk{{\hat\tau}_{m,k}}
\def\htaumkj{{\hat\tau}_{m,k,j}}
\def \IJK{\mathcal I}
\def \bga{\bar\gamma}
\def\taup{{\oo{\tau}}}
\def\sip{{\oo{\zeta}}}
\def\sii{{\zeta}}

\def\P{{\Bbb P}}
\def\E{{\Bbb E}}
\def\NNN{{\mathcal N}}
\def \njc{{\bf !!! note change  !!! }}
\def\ol#1{{\overline{#1}}}

\bibliographystyle{amsplain}

\title[Frequent Points for Random Walks] {Frequent Points for Random Walks
in Two Dimensions}

\author[Richard F. Bass\quad \mbox{\rm and}\quad Jay Rosen] {Richard F.
Bass$^*$\quad \mbox{\rm and }
         Jay Rosen$^\dagger$}

\date{July 25, 2006.
\newline\indent
$^*$Research  partially supported by NSF grant \#DMS-0601783.
\newline\indent
$^\dagger$Research supported, in part, by grants from the NSF and from
PSC-CUNY}

\begin{abstract}
\noindent For a symmetric random walk in
$\Z^2$ which does not necessarily have bounded  jumps
         we study those points which are
         visited an unusually large number of times. We prove the 
analogue of the
Erd\H{o}s-Taylor  conjecture and obtain the asymptotics for the number of
visits to the most visited site. We also obtain the asymptotics for 
the number of
points which are visited very frequently by time $n$.
   Among the tools we use are Harnack inequalities and Green's function
estimates for random walks with unbounded jumps; some of these are of
independent interest.
\end{abstract}

\maketitle

\section{Introduction}

The paper \cite{DPRZ4} proved a conjecture of Erd\H{o}s and  Taylor
concerning the number   $L^{ \ast}_n$ of visits to the most visited site for
simple random walk in
$\Z^2$ up to step $n$. It was  shown there that
\begin{equation}
\lim_{n\to \ff}{L^{ \ast}_n\over (\log n)^2}=1/\pi
\hspace{.1in}\mbox{a.s.}\label{0.0}
\end{equation}
         The approach in that paper was to first prove an analogous result for
planar  Brownian motion and then to use strong approximation. This approach
applies  to other random walks, but only if they have moments of all orders. In
a more recent paper
\cite{jhung}, Rosen developed purely random walk methods which allowed him
to prove (\ref{0.0}) for simple random walk. A key to the approach both for
Brownian motion and simple random walk is to exploit a certain tree structure
with regard to excursions between nested families of disks. When we turn to
random walks with jumps, this tree structure  is no longer automatic, since the
walk may jump across disks. In this paper we show how to extend the method  of
\cite{jhung} to symmetric recurrent  random walks
$X_j,\,j\geq 0$, in
$\Z^2$. Not surprisingly, our key task is to control the jumps across 
disks. Our
main conclusion is that it suffices to require that for some
$\bb>0$
\begin{equation}
\E|X_1|^{3+2\beta}<\infty,\label{1.0}
\end{equation} together with some mild uniformity conditions. (It will make
some formulas later on  look nicer if we use
$2\bb$ instead of
$\bb$ here.) We go beyond  (\ref{0.0}) and study the size of the set of
`frequent points,' i.e. those points in
$\Z^2$ which are
         visited an unusually large number of times,
         of order $(\log n)^2$. Perhaps more important than our 
specific results,
we develop powerful estimates for our random walks which we expect will have
wide applicability. In particular, we develop Harnack inequalities extending
those of
\cite{L3} and we develop estimates for Green's functions for random walks
killed on entering a disk. The latter estimates are new even for simple random
walk and are of independent interest.

We assume for simplicity that $X_{ 1}$ has
         covariance matrix equal to the identity and that $X_n$ is strongly
aperiodic. Set
$p_{ 1}(x,y)=p_{ 1}(x-y)=\P^{ x}(X_1=y )$. We will say that our walk satisfies
Condition A if the following holds.
\medskip

\noindent{\bf Condition A.} {\sl Either $X_{ 1}$ is  finite range, that is,
$p_{ 1}( x)$ has bounded support, or else for any  $s\leq R$ with $s $
sufficiently large
\begin{equation}
\inf_{y\,;\,R\leq |y|\leq R+s}\,\,
\sum_{ z\in D(0,R)}p_{ 1}(y,z)\geq ce^{ -\bb\,s^{1/4}}.\label{condA}
\end{equation} } Condition A is implied by
\begin{equation} p_{ 1}(x)\geq ce^{ -\bb\,|x|^{1/4}},\hspace{ .2in}x\in
\Z^2,\label{condA.3}
\end{equation} but  (\ref{condA}) is much weaker. (\ref{condA}) is a mild
uniformity condition, and is used to obtain Harnack inequalities. Recent work
on Harnack inequalities for  processes with jumps indicates that some sort of
uniformity condition is needed; see
\cite{BK}.

          Let $L_n^{ x}$ denote the number of times that $x\in \Z^2$ 
is visited by
         the random walk in $\Z^2$ up to step $n$ and set $L^{
\ast}_n=\max_{ x\in
\Z^2}L_n^{ x}$.

\begin{theorem}\label{theo-et} Let $\{ X_{ j}\,;\,j\geq 1\}$ be a symmetric
strongly aperiodic  random walk in
$\Z^2$  with $X_1$ having the identity as the covariance matrix and satisfying
Condition A and (\ref{1.0}). Then
\begin{equation}
\lim_{n\to \ff}{L^{ \ast}_n\over (\log n)^2}=1/\pi
\hspace{.1in}\mbox{a.s.}\label{0.0h}
\end{equation}
\end{theorem}

Theorem \ref{theo-et} is the analogue of the Erd\H{o}s-Taylor conjecture for
simple random walks. We also look at how many frequent points there are.
         Set
\begin{equation}
\Theta_n(\al)=\Big\{x\in  \Z^2:\;
         \frac{L^{ x}_{ n}}{(\log n)^2}\geq \al/\pi
\Big\}.\label{bthick.1a}
\end{equation}

For any set
$B\subseteq \Z^2$ let
$T_{B}=\inf\{i\geq 0\,|\,X_i\in B\}$ and let
$|B|$ be the cardinality of $B$. Let
\begin{equation}
\Psi_n(a)=\Big\{x\in  D(0, n):\;
         \frac{L^{ x}_{ T_{D(0, n)^{ c}}}}{(\log n)^2}\geq 2a/\pi
\Big\}\label{bthick.1}
\end{equation}

\begin{theorem}
\label{theo-late} Let $\{ X_{ j}\,;\,j\geq 1\}$ be as in Theorem
\ref{theo-et}.  Then for any
$0<\al< 1$
\begin{equation}
\lim_{n\to \ff}{\log |\Theta_n(\al)|\over \log n}=1-\al
\hspace{.1in}\mbox{a.s.}\label{thick.2}
\end{equation} Equivalently, for any $0< a< 2$
\begin{equation}
\lim_{n\to \ff}{\log |\Psi_n(a)|\over \log n}=2-a
\hspace{.1in}\mbox{a.s.}\label{bthick.2}
\end{equation}
\end{theorem}

The equivalence of (\ref{thick.2}) and (\ref{bthick.2}) follows from 
the fact that
\begin{equation}
         \lim_{n\to \ff}\frac{\log T_{D(0, n)^{ c}}}{\log n}=2
\hspace{.1in}\mbox{a.s.}\label{equiv.1z}
\end{equation} For the convenience of the reader we give a proof of 
this fact in
the appendix.

In Section \ref{sec-rwprelim} we collect some facts about random walks in
$\Z^2$, and in Section \ref{sec-harn} we establish the Harnack inequalities we
       need. The upper bound of Theorem \ref{theo-late} is proved in Section
       \ref{sec-hit}. The lower bound is established in Section
\ref{sec-lowerbound}, subject to certain estimates which form the 
subject of the
following three sections. An appendix gives further information about random
walks in $\Z^2$.

   There is a good deal of flexibility in our choice of Condition A. 
For example, if
$\E|X_1|^{4+2\beta}<\ff$, we can replace $\bb\,s^{1/4}$ by
$s^{1/2}$. On the other hand, if we merely assume  that
$\E|X_1|^{2+2\beta}<\ff$, our methods do not allow us to obtain  any useful
analogue of the Harnack inequalities we derive in Section
\ref{sec-harn}.

\section{Random Walk Preliminaries}\label{sec-rwprelim}

Let $X_n,\,n\geq 0$, denote a symmetric recurrent random walk in
$\Z^2$ with covariance matrix equal to the identity.  We set
$p_{ n}(x,y)=p_{ n}(x-y)=\P^{ x}(X_n=y )$ and assume that for some
     $\bb>0$
\begin{equation} \E|X_1|^{3+2\beta}=\sum_x |x|^{3+2\beta} p_1(x)<\infty.
\label{2.0}
\end{equation} (It will make some formulas later on  look nicer if we use
$2\bb$ instead of
$\bb$ here.)
         In this section we collect some facts about
$X_n,\,n\geq 0$, which will be used in our paper. The estimates for 
the interior
of a ball are the analogues for 2 dimensions and
$3+2\beta$ moments of some results that are proved in \cite{L3} for random
walks in dimensions 3 and larger that have $4$ moments. Several results which
are well known to experts but are not in the literature are given in 
an appendix.

We will assume throughout that $X_n$ is strongly aperiodic.
         Define the potential kernel
\begin{equation} a(x)=\lim_{ n\rar\ff}\sum_{ j=0}^{ n}\lc p_{ j}(0)-p_{
j}(x)\rc.\label{p2.1}
\end{equation} We have that $a(x)<\ff$ for all $x\in\Z^2$,
$a(0)=0$, and for
$|x|$ large
\begin{equation} a(x)={ 2\over \pi}\log |x| +k +o( |x|^{ -1})\label{p2.2}
\end{equation} with $k$ an explicit constant. See Proposition \ref{AP2} for a
proof.  We also note that since $p_{ j}(x)\leq p_{ j}(0)$ for any symmetric
random walk, we have $a(x)\geq 0$.

Let $D(x,r)=\{y\in\Z^2\,|\,|y-x|<r\} $. For any set $A\subseteq
\Z^2$  we define the boundary $\partial A$ of $A$ by $\partial
A=\{y\in\Z^2\,|\,y\in A^c,\,\mbox{and }
\inf_{x\in A}|y-x|\leq 1\} $ and the
$s$-band $\partial A_{ s}$ of  $\partial A$ by
$\partial A_{ s}=\{y\in\Z^2\,|\,y\in A^c,\,\mbox{and } \inf_{x\in A}|y-x|\leq
s\} $. For any set
$B\subseteq \Z^2$ let
$T_{B}=\inf\{i\geq 0\,|\,X_i\in B\}$. For
$x,y\in A$ define the Green's function
\begin{equation} G_A(x,y)=\sum_{i=0}^\ff
\E^x\(X_i=y,\,i<T_{A^c}\).\label{p1.1}
\end{equation}

For some $c<\ff$
\begin{equation}
\E^x(T_{D(0,n)^c})\leq c n^2, \qquad x\in D(0,n),\hspace{ .2in} n\geq
1.\label{Lesc.1}
\end{equation} This is proved in Lemma \ref{ALesc}. In particular,
\begin{equation}
\sum_{y\in D(0,n)} G_{ D(0,n)}(x,y)\leq cn^2.\label{Lesc.2}
\end{equation}

       Define the hitting distribution of $A$ by
\begin{equation} H_{A}(x,y)=\P^x(X_{T_{A}}=y).\label{p1.4}
\end{equation}
       As in Proposition 1.6.3 of \cite{L1}, with $A$ finite, by 
considering the first
hitting time of $A^c$ we have that  for $x,z\in A$
\begin{equation} G_A(x,z)=\lc \sum_{y\in A^c}H_{A^c }(x,y)a( y-z)\rc-a(
x-z).\label{p2.3}
\end{equation} In particular
\begin{eqnarray} G_{D(0,n) }(0,0)&= &\sum_{n\leq |y|\leq n+n^{ 3/4}}H_{
D(0,n)^{ c} }(0,y)a( y)
           \label{p2.4}\\
         &&+\sum_{|y|> n+n^{ 3/4}}H_{D(0,n)^{ c} }(0,y)a( y).
           \nn
\end{eqnarray} Using the last exit decomposition
\begin{equation} H_{D(0,n)^{ c} }(0,y)=\sum_{ z\in D(0,n)}G_{ D(0,n)}(0,z)p_{
1}( z,y)
\label{p2.4m1}
\end{equation} together with (\ref{Lesc.2}) and (\ref{2.0}) we have for any
$k\geq 1$
\begin{eqnarray}
\sum_{|y|\geq n+kn^{ 3/4} }H_{D(0,n)^{ c} }(0,y)&\leq &
\sum_{ z\in D(0,n)}G_{ D(0,n)}(0,z)P( |X_{ 1}|\geq kn^{ 3/4})\nn\\ &\leq &
      cn^{ 2}/( kn^{ 3/4})^{3+2\bb}
\leq c/(k^{ 3}n^{1/4+\bb}).\label{p2.4m}
\end{eqnarray}

Using this together with (\ref{p2.2})  we can bound the last term in 
(\ref{p2.4})
by
\begin{eqnarray} &&
\sum_{|y|> n+n^{ 3/4}}H_{D(0,n)^{ c} }(0,y)a(y)
           \label{p2.4a}\\
         &&\hspace{.2in}
\leq\sum_{ k=1}^{ \ff}{ c\over k^{ 3}n^{1/4+\bb}}\sup_{n+ kn^{ 3/4}\leq |x|
\leq n+(k+1)n^{3/4}}a(x)\nn\\
         &&\hspace{.2in}
\leq { C\over n^{1/4+\bb}}\sum_{ k=1}^{ \ff}{ 1\over k^{ 3}}\log
(n+(k+1)n^{3/4})=O(n^{ -1/4}).\nn
\end{eqnarray}
       Using (\ref{p2.2}) for the first term in (\ref{p2.4}) then gives
\begin{equation} G_{D(0,n) }(0,0)={ 2\over \pi}\log n +k+O(n^{
-1/4}).\label{p2.5}
\end{equation}

Let $\eta=\inf\{i\geq 1\,|\,X_i\in \{  0 \}\cup D(0,n)^c\}$. Applying the
optional sampling theorem to the martingale $a( X_{ j\wedge
\eta})$ and letting $j\to \infty$, we have that for any
$x\in D(0,n)$
\begin{equation} a( x)=\E^{ x}(a( X_{\eta}))=\E^{ x}(a(
X_{\eta})\,;\,X_{\eta}\in D(0,n)^c).\label{p2.6}
\end{equation}   To justify taking the limit as $j\to \infty$, note that
    $|a(X_{j\land \eta})|^2$ is a submartingale, so
$\E|a(X_{j\land\eta})|^2\leq \E|a(X_\eta)|^2$, which is finite by (\ref{p2.2})
and (\ref{p2.4m}); hence the family of random variables
$\{a(X_{j\land
\eta})\}$ is uniformly integrable. Using (\ref{p2.2}) and the analysis of
(\ref{p2.4a})  we find that
\begin{eqnarray} && \E^{ x}\(a( X_{\eta})\,;\,X_{\eta}\in D(0,n)^c\)
           \label{ p2.7}\\
         &&=\sum_{y\in \partial D(0,n)_{n^{ 3/4}}} a(
y)\PPP^x\(X_{\eta}=y\)+\sum_{y\in D(0,n+n^{ 3/4})^c} a(
y)\PPP^x\(X_{\eta}=y\)
           \nn\\
         &&=\({ 2\over \pi}\log n +k +O(n^{ -1/4})\)\P^{ x}\(\,X_{\eta}\in
D(0,n)^c\)+O(n^{ -1/4}).
           \nn
\end{eqnarray}  Using this and (\ref{p2.2})  we find that for
$0<|x|<n$,
\begin{eqnarray}
         \PPP^x\(T_{0}<T_{D(0,n)^c}\) &=& {(2/ \pi)\log (n/ |x|)+O( |x|^{-1/4})
\over (2/ \pi)\log n +k +O(n^{ -1/4})} \nn\\
          &=& {\log (n/ |x|)+O( |x|^{-1/4} )\over
\log (n)}( 1+O( (\log n)^{-1}).  \label{p1.2c}
\end{eqnarray}

By the strong Markov property
\begin{equation} G_{D(0,n)}(x,0)=\PPP^x\(T_{0}<T_{D(0,n)^c}\)
G_{D(0,n)}(0,0).\label{mark.1}
\end{equation}
    Using this, (\ref{p2.5}) and the first line of (\ref{p1.2c}) we obtain
\begin{equation} G_{D(0,n)}(x,0)={2\over \pi}\log \Big({n\over |x|}\Big) +O(
|x|^{-1/4} ).
\label{p1.2}
\end{equation} Hence
\begin{equation} G_{D(0,n)}(x,y)\leq G_{D(x,2n)}(x,y)\leq c\log n.
\label{p1.2j}
\end{equation}

   Let $r<|x|<R$ and $\zeta=T_{D(0,R)^c}\land T_{D(0,r)}$. Applying the
argument leading to   (\ref{p1.2c}), but with the martingale $a(X_{j\land
\zeta})-k$  we can obtain that  uniformly in $r<|x|<R$
\begin{equation}
\PPP^x\(T_{D(0,R)^c}<T_{ D(0,r)}\)= {\log (|x|/r)+O( r^{-1/4})\over
\log (R/ r)}\label{p1.2bfa}
\end{equation} and
\begin{equation}
\PPP^x\(T_{ D(0,r)}<T_{D(0,R)^c}\)= {\log (R/ |x|)+O( r^{-1/4})\over
\log (R/ r)}.\label{p1.2b}
\end{equation}

Using a last exit decomposition, for any $0<\delta<\epsilon<1$, uniformly in
$x\in D(0,n)\setminus D(0,\epsilon n)$
\begin{align}
\P^x(|X_{T_{D(0,\epsilon n)}\land T_{D(0,n)^c}}|\leq \delta n) &=\sum_{z\in
D(0,n)\setminus D(0,\epsilon n)}
\sum _{w\in D(0,\delta n)} G_{D(0,n)}(x,z) p_1(z,w)\nn\\ &\leq cn^2
\log n
\P(|X_1|>(\epsilon-\delta) n)\nn\\ &\leq cn^2 \log
n\frac{1}{n^{3+2\beta}}=c\log n\frac{1}{n^{1+2\beta}}.\label{3m.2}
\end{align} Here we used (\ref{p1.2j}) and the fact that
     $|z-w|\geq (\epsilon-\delta) n$.

   We need a more precise error term when $x$ is near the outer boundary.
   Let $\rho(x)=n-|x|$. We have the following lemma.

\begin{lemma}\label{lem-gamblow} For any $0<\delta<\epsilon<1$ we can
find $0<c_{ 1}<c_{ 2}<\infty$, such that for all
$x\in D(0, n)\setminus D(0, \epsilon n)$ and all $n$ sufficiently large
\begin{equation} c_{ 1}\frac{\rho( x)\lor 1}{n}\leq \P^{ x}(T_{D(0, \delta n)
}<T_{D(0, n)^{ c} })\leq c_{ 2}\frac{\rho( x)\lor 1}{n} .\label{92.1}
\end{equation}
\end{lemma}

{\bf  Proof: } Upper bound: By looking at two lines, one tangential to
$\partial D(0,n)$ and perpendicular to the ray from 0 to $x$, and the other
parallel to the first but at a distance $\delta n$ from 0, the upper 
bound follows
from the gambler's ruin estimates of \cite[Lemma 2.1]{L3}.

Lower bound: We first show that for any $\zeta>0$ we can find a constant $c_{
\ze}>0$ such that
\begin{equation} c_{ \ze}\frac{\rho( x)}{n}\leq \P^{ x}(T_{D(0, \delta n)
}<T_{D(0, n)^{ c} })\label{92.1z}
\end{equation}
   for all $  x\in D(0, n-\ze)\setminus D(0, \epsilon n)$.

   Let $T=\inf\{ t\,|\,X_{ t}\in D(0, \delta n)\cup D(0,n)^{ c}\}$ and
$\gamma\in (0,\frac15)$. Let $\ol a(x)=\frac{2}{\pi}(a(x)-k)$, where $k$ is the
constant in (\ref{p2.2}) so that $\ol a(x)=\log |x|+o(1/|x|)$.  Clearly
$\ol a(X_{j\land T})$ is a martingale, and by
      the optional sampling theorem,
\begin{eqnarray} && \ol a( x)=\E^{ x}( \ol a( X_{ T});\,X_{ T}\in D(0, n)^{
c})\label{92.2}\\ &&
\hspace{ .4in}+\E^{ x}( \ol a( X_{ T});\,X_{ T}\in D(0, \delta n)\setminus D(0,
\delta n/2))\nn\\ &&
\hspace{ .4in} +\E^{ x}( \ol a( X_{ T});\,X_{ T}\in D(0, \delta n/2)).\nn
\end{eqnarray}  It follows from  (\ref{3m.2}) that
\begin{equation} \E^{ x}( \ol a( X_{ T});\,X_{ T}\in D(0, \delta n/2))=O(n^{
-1-\gamma})\label{92.5}
\end{equation}  and
\begin{equation}\label{CII1}
\P^x(X_T\in D(0,\delta n/2))=O(n^{-1-\gamma}).
\end{equation} From (\ref{92.2}) we see that
\begin{align}
\log |x|&\geq \log n\,\P^x(T_{D(0,n)^c}<T_{D(0, \delta n)})\label{92.6}\\
&~~~+\log(\delta n/2) \P^x( T_{D(0,n)^c}>T_{D(0, \delta n)})+o(1/n)\nn\\
&= \log n [1-\P^x(T_{D(0,\delta n)}< T_{D(0,  n)^c})] +
\log(\delta n/2)
\P^x(T_{D(0,\delta n)}< T_{D(0,  n)^c})]\nn\\ &~~~~ +o(1/n)\nn\\ &=\log n
+(\log(\delta n/2)-\log n)
\P^x(T_{D(0,\delta n)}< T_{D(0,  n)^c})+o(1/n).\nn
\end{align}
 Note that for some $c>0$
\begin{equation}
\log (1-z)\leq -cz,\hspace{ .2in}0\leq z\leq 1-\ep.\label{92.patch}
\end{equation}
Hence for $x\in D(0, n-\ze)\setminus D(0,\epsilon n)$
$$\log(n/|x|)=-\log\Big(1-\frac{(n-|x|)}{n}\Big)\geq c\frac{\rho( x)}{n}$$
  Solving (\ref{92.6}) for
$\P^x(T_{D(0,\delta n)}< T_{D(0,  n)^c})$ and using the fact that
and $\rho(x)\geq \ze$ to control the $o(1/n)$ term completes the proof of
(\ref{92.1z}).

Let $A=D(0, n-\ze)\setminus D(0,\epsilon n)$. Then by the strong Markov
property and (\ref{92.1z}), for any $x\in D(0,n)\setminus D(0,n-\ze)$
\begin{eqnarray} & &\P^{ x}(T_{D(0, \delta n) }<T_{D(0, n)^{ c}
})\label{92.20}\\  &&\geq \P^{ x}(T_{D(0, \delta n) }\circ\theta_{
T_{A}}<T_{D(0, n)^{ c} }\circ\theta_{ T_{A}};\,  T_{A}<T_{D(0, n)^{ c} })
\nonumber\\  &&=\E^{ x}(\P^{X_{ T_{A}}}(T_{D(0, \delta n) }<T_{D(0, n)^{ c}
});\,  T_{A}<T_{D(0, n)^{ c} })
\nn\\  &&\geq c_{ \ze}\,\,{\ze \over n}\,\, \P^{ x}( T_{A}<T_{D(0, n)^{ c}}).
\nn
\end{eqnarray}
   (\ref{92.1}) then follows if we can find $\ze>0$ such that uniformly in $n$
\begin{equation}
\inf_{ x\in D(0,n)\setminus D(0,n-\ze)} \P^{ x}( T_{A}<T_{D(0, n)^{
c}})>0.\label{92.21}
\end{equation} To prove (\ref{92.21}) we will show that we can find
$N<\infty$ and
$\zeta, c>0$ such that for any $x\in\Z^{ 2}$ with $|x|$ sufficiently 
large we can
find $y\in\Z^{ 2}$ with
\begin{equation} p_{ 1}( x,y)\geq c\hspace{ .2in}\mbox{ and }\hspace{
.2in}|x|-N\leq |y|\leq |x|-\zeta.\label{92.30}
\end{equation}

Let $C_{ x}$ be the cone with vertex at $x$ that contains the origin, has
aperture $9\pi/10$ radians, and such that the line through 0 and $x$ bisects
the cone. If $|x|$ is large enough,
$C_{ x}\cap D(x,N)$ will be contained in $D(0,|x|)$. We will show 
that there is a
point $y\in \Z^2\cap C_{ x}\cap D(x ,N)$, $y\neq x$, which satisfies
(\ref{92.30}). Note that for any $y\in \Z^2\cap C_{ x}\cap D(x ,N)$,
$y\neq x$, we have
\begin{equation} 1\leq |y-x|\leq N.\label{92.31}
\end{equation} Furthermore, if $\al$ denotes the angle between the line from
$x$ to the origin and the line from $x$ to $y$, by the law of cosines,
\begin{equation} |y|^{ 2}=|x|^{ 2}+|y-x|^{ 2}-2\cos (\al)\,|x|\, |y-x|.
\label{92.32}
\end{equation} Then for $|x|$ sufficiently large, using (\ref{92.31}),
\begin{eqnarray}  |y|&=&\sqrt{|x|^{ 2}+|y-x|^{ 2}-2\cos (\al)\,|x| \,|y-x|}
\label{92.33}\\ &=& |x|\sqrt{1+{ |y-x|^{ 2}\over |x|^{ 2}}-{2\cos (\al)\,|y-x|
\over |x|}}\nonumber\\ &=& |x|
\left(1-{\cos (\al)\,|y-x| \over |x|}+O( { 1\over |x|^{ 2}})\right)\nonumber\\
     &=& |x| -\cos (\al)\,|y-x|+O( { 1\over |x|})\nonumber\\
     &\leq & |x| -\cos (9\pi/20)/2.\nonumber
\end{eqnarray} Setting $\zeta=\cos (9\pi/20)/2>0$ we see that the second
condition in (\ref{92.30}) is satisfied for all $y\in \Z^2\cap C_{ x}\cap D(x
,N)$, $y\neq x$, and it suffices to show that we can find such a $y$ with
$p_{ 1}( x,y)\geq c$. ($c,N$ remain to be chosen).

By translating by $-x$ it suffices to work with cones having their 
vertex at the
origin. We let $C( \theta,\theta')$ denote the cone with vertex at the origin
whose sides are the half lines making angles
$\theta<\theta'$  with the positive $x$-axis. Set
$C( \theta,\theta',N)=C( \theta,\theta')\cap D(0 ,N)$. It suffices to show that
for any $ \theta$ we can find $y\in  C(
\theta,\theta+9\pi/10,N)$, $y\neq 0$, with $p_{ 1}( 0,y)\geq c$. Let
$j_{ \th}=\inf\{ j\geq 0\,|\,j\pi/5 \geq \th\}$. Then it is easy to see that
\[C(j_{ \th}\pi/5,j_{ \th}\pi/5+2\pi/3,N)\subseteq C(
\theta,\theta+9\pi/10,N).\] It now suffices to show that
     for each $0\leq j\leq 9$ we can find
$y_{ j}\in C(j\pi/5,j\pi/5+2\pi/3)$, $y_{ j}\neq 0$, with $p_{ 1}( 
0,y_{ j})>0$,
since we can then take
\begin{equation} c=\inf_{ j }p_{ 1}( 0,y_{ j})\hspace{ .2in}\mbox{ and
}\hspace{ .2in} N=2\sup_{ j}|y_{ j}|.\label{92.34}
\end{equation}

First consider the cone $C(-\pi/3,\pi/3 )$, and recall our assumption that the
covariance matrix of $X_{ 1}=(X_{ 1}^{(1)},X_{ 1}^{(2)} )$ is $I$. If
\[\P(X_{ 1}\in C(-\pi/3,\pi/3 ),\,X_{ 1}\neq 0)=0\] then by symmetry,
$\P(X_{ 1}\in -C(-\pi/3,\pi/3 ),\,X_{ 1}\neq 0)=0$. Therefore \[
\P\Big(|X_{ 1}^{(2)}| > |X_{ 1}^{(1)}|\,\Big|\,X_{ 1}\neq 0\Big)=1.\] But then
$1=\E |X_{ 1}^{(2)}|^2> \E|X_{ 1}^{(1)}|^2=1$, a contradiction. So there must
be  a point  $y\in C(-\pi/3,\pi/3 )$, $y\neq 0$, with $p_{ 1}( 0,y)>0$.

Let
$A_j$ be the rotation matrix such that the image of
$C(j\pi/5,j\pi/5+2\pi/3)$ under
$A_j$ is the cone $C(-\pi/3,\pi/3 )$ and let
     $Y_{ 1}=A_jX_{ 1}$. Then
\begin{eqnarray}&&
\P(X_1\in C(j\pi/5,j\pi/5+2\pi/3),\,X_{ 1}\neq 0)\label{92.35}\\ &&=\P(Y\in
C(-\pi/3,\pi/3 ),\,Y_{ 1}\neq 0).\nn
\end{eqnarray}
      Note
$Y_{ 1}$ is mean 0, symmetric, and since $A_j$ is a rotation matrix, the
covariance matrix of
$Y_{ 1}$ is the identity. Hence the argument of the last paragraph shows that
the probability in (\ref{92.35}) is non-zero. This completes the proof of
  of our lemma.
     \qed

\begin{lemma}\label{lem-greblow2} For any $0<\delta<\epsilon<1$ we can
find $0<c_{ 1}<c_{ 2}<\infty$, such that for  all $x\in D(0, n)\setminus D(0,
\ep n), y\in D(0, \de n)$ and all $n$ sufficiently large
\begin{equation} c_{ 1}\frac{\rho( x)\lor 1}{n}\leq G_{D(0, n) }( y,x)\leq c_{
2}\frac{\rho( x)\lor 1}{n} .\label{92.13}
\end{equation}
\end{lemma}

{\bf  Proof:} Upper bound: Choose $\delta<\ga<\epsilon$ and let
$T'=\inf\{ t\,|\,X_{ t}\in D(0, \ga n)\cup D(0,n)^{ c} \}$. By the 
strong Markov
property,
\begin{eqnarray} \hspace{ .2in} G_{D(0, n) }(x,y)&=&\sum_{z\in  D(0,
\ga n)}\P^{ x}(X_{ T'}=z )G_{D(0, n) }(z,y)\nn\\ &=&\sum_{z\in  D(0,
\de n)}\P^{ x}(X_{ T'}=z )G_{D(0,n) }(z,y)\nn\\ &+&\sum_{z\in  D(0,
\ga n)\setminus D(0, \de n)}\P^{ x}(X_{ T'}=z )G_{D(0,n) }(z,y)\nn\\ &\leq
&c\P^{ x}(X_{ T'}\in D(0, \de n) )\log n+c\P^{ x}(X_{ T'}\in D(0,
\ga n) ).\label{92.15}
\end{eqnarray} Here we used (\ref{p1.2j}) and the bound
$G_{D(0,n) }(z,y)\leq G_{D(y,2n) }(z,y)\leq c$ uniformly in $n$ which follows
from (\ref{p1.2}) and the fact that
$|z-y|\geq ( \ga-\de)n$. The upper bound in (\ref{92.14}) then follows from
(\ref{3m.2}) and  Lemma
\ref{lem-gamblow}.

Lower bound:  Let
$T=\inf\{ t\,|\,X_{ t}\in D(0, \de n)\cup D(0,n)^{ c} \}$. By the strong Markov
property,
\begin{equation} G_{D(0, n) }(x,y)=\sum_{z\in  D(0, \de n)}\P^{ x}(X_{ T}=z
)G_{D(0, n) }(z,y).\label{92.14}
\end{equation} The lower bound in (\ref{92.14}) follows  from Lemma
\ref{lem-gamblow} once we show that
\begin{equation}
\inf_{y,z\in  D(0, \de n)}G_{D(0, n)}(z,y)\geq a>0\label{92.14r}
\end{equation} for some $a>0$ independent of
$n$. We first note that
\begin{eqnarray} &&
\inf_{y,z\in  D(0, \de n)}\inf_{z\in  D(y, (\ep-\de) n/2)}G_{D(0, n)}(z,y)
\label{92.14s}\\ &&\hspace{ .3in}\geq \inf_{y\in  D(0, \de n)}\inf_{z\in D(y,
(\ep-\de)  n/2)}G_{D(y, (\ep-\de) n)}(z,y)>0\nn
\end{eqnarray}
     uniformly in $n$  by (\ref{p1.2}). But by the invariance 
principle, there is a
positive probability independent of
$on$ that the random walk starting at $y$ will enter  $D(y, (\ep-\de) n/2)\cap
D(0, \de n)$ before exiting $D(0, n)$. (\ref{92.14r}) then follows from
(\ref{92.14s}) and the strong Markov property.
\qed

\begin{lemma}\label{lem-jumpo}
\begin{equation}
\sup_{x\in D(0,n/2)}\PPP^{x}(T_{ D(0,n)^c}
\neq T_{\partial D(0,n)_{ s}})\leq c(s^{-1-2\bb}\vee n^{-1-2\bb}\log
n).\label{c2.11k}
\end{equation}
\end{lemma}
\noindent {\bf Proof of Lemma~\ref{lem-jumpo}:} We begin  with the last exit
decomposition
\begin{eqnarray}
         &&\sup_{x\in D(0,n/2)}\PPP^{x}(T_{ D(0,n)^c}
\neq T_{\partial D(0,n)_{ s}})
           \label{jumpo.1}\\
         &&=\sup_{x\in D(0,n/2) }
\sum_{\stackrel{y\in D(0,n)}{w\in  D(0,n+s)^{ c}}} G_{D(0,n) }(x,y 
)p_{ 1}(y,w )
           \nn\\
         &&=\sup_{x\in D(0,n/2) }
\sum_{\stackrel{|y|\leq 3n/4}{n+s\leq |w|}} G_{D(0,n) }(x,y )p_{ 1}(y,w )
           \nn\\
         &&+\sup_{x\in D(0,n/2) }
\sum_{\stackrel{3n/4<|y|< n}{n+s\leq |w|}} G_{D(0,n) }(x,y )p_{ 1}(y,w ).
           \nn
\end{eqnarray} Using (\ref{p1.2j}) and (\ref{2.0})
\begin{eqnarray}
         &&\sup_{x\in D(0,n/2) }
\sum_{\stackrel{|y|\leq 3n/4}{n+s\leq |w|}} G_{D(0,n) }(x,y )p_{ 1}(y,w )
           \label{jumpo.2}\\
         &&\leq c \log n
\sum_{ |y|\leq 3n/4 }
\P(|X_{ 1}|\geq n/4 )
           \nn\\
         &&\leq c \log n
\sum_{|y|\leq 3n/4} { 1\over |n |^{ 3+2\bb}}
\leq cn^{-1-2\bb}\log n.
           \nn
\end{eqnarray} Using (\ref{92.13}) and (\ref{2.0}) we
     have
\begin{eqnarray}
         &&\sup_{x\in D(0,n/2) }
\sum_{\stackrel{3n/4<|y|< n}{n+s<|w|}} G_{D(0,n) }(x,y )p_{ 1}(y,w )
           \label{jumpo.3}\\
         &&\leq c n^{-1}
\sum_{3n/4<|y|\leq n} ( n-|y|)\P(|X_{ 1}|\geq s+n-|y| )
           \nn\\
         &&\leq c n^{-1}
\sum_{3n/4<|y|\leq n} { n-|y|\over (s+n-|y|)^{ 3+2\bb}}
           \nn\\
         &&
\leq c n^{-1}
\sum_{3n/4<|y|< n} { 1\over (s+n-|y|)^{ 2+2\bb}}\leq cs^{-1-2\bb}.
           \nn
\end{eqnarray} Here, we bounded the last sum by an integral and used polar
coordinates:
\begin{eqnarray}\qquad
\sum_{3n/4<|y|< n} { 1\over (s+n-|y|)^{ 2+2\bb}}&\leq  & 2\pi\int_{3n/4 }^{
n}{ 1\over (s+n-u)^{ 2+2\bb}}\,u\,du\label{jumpo.4}\\ &\leq  & cn\int_{0}^{
n/4}{ 1\over (s+u)^{ 2+2\bb}}\,du\leq cns^{ -1-2\bb}.
\nonumber
\end{eqnarray}
\qed

Combining (\ref{p1.2c}) and (\ref{c2.11k}) we obtain, provided
$n^{-1-2\bb}\log n\leq s^{-1-2\bb}$ and $|x|\leq n/2$,
\begin{eqnarray}
         &&
         \PPP^x\(T_{D(0,n)^c}<T_{0}\,;\,T_{ D(0,n)^c} =T_{\partial 
D(0,n)_{ s}}\)
\label{p1.2bb}\\
         &&
         =1-{(2/ \pi)\log (n/ |x|)+O( |x|^{-1/4})
\over (2/ \pi)\log n +k +O(n^{ -1/4})}+O(s^{-1-2\bb}).
\nn
\end{eqnarray}

     Using (\ref{p1.2bfa}) and (\ref{c2.11k}) we then obtain, provided
$R^{-1-2\bb}\log R\leq s^{-1-2\bb}$ and $|x|\leq R/2$,
\begin{eqnarray}
         &&
         \PPP^x\(T_{D(0,R)^c}<T_{ D(0,r)}\,;\,T_{ D(0,R)^c} 
=T_{\partial D(0,R)_{
s}}\)
\label{p1.2bd}\\
         &&
         = {\log (|x|/r)+O( r^{-1/4})\over
\log (R/ r)}+O(s^{-1-2\bb}).  \nn
\end{eqnarray}

\begin{lemma}\label{lem-extjumpsmall}For any $s<r<R$ sufficiently large with
$ R\leq  r^{ 2}$ we can find
$c<\infty$ and
$\delta>0$ such that for any $r<|x|<R$
\begin{equation}
\P^x(T_{D(0,r)}<T_{D(0,R)^c}; X_{T_{D(0,r)}}\in  D(0,r-s))
\leq cr^{ -\delta}+cs^{-1-2\beta}.
\label{G10a}
\end{equation}
\end{lemma}

\proof
Let $A( R,r)$ denote the annulus $D(0,R)\setminus D(0,r)$.
Using a last exit decomposition  we have
\begin{align}
\P^x(T_{D(0,r)}<T_{D(0,R)^c}; X_{T_{D(0,r)}}&\in  D(0,r-s))
\nn\\
&=\sum_{w\in D(0,r-s)}~~\sum_{y\in A} G_A(x,y)p_1(y,w)
\nn\\
&=\sum_{w\in D(0,r-s)}~~\sum_{r< |y|\leq r+r^{1/(2+\bb)}} G_A(x,y)p_1(y,w)
\nn\\
&+\sum_{w\in D(0,r-s)}~~\sum_{r+r^{1/(2+\bb)}< |y|<R} G_A(x,y)p_1(y,w).
\label{G10}
\end{align}
By (\ref{2.0}), for $y\in A$
$$\sum_{w\in D(0,r-s)} p_1(y,w)\leq \frac{c}{(|y|-(r-s))^{3+2\beta}}.$$
Let  $U_k=\{y\in \Z^2:
r+k-1< |y|\leq r+k\}$. We show below that we can find $c<\ff$ such that
for all  $1\leq k\leq r^{1/(2+\bb)}$,
\begin{equation}
\label{G6}
\sum_{y\in U_k} G_A(x,y)\leq ck, \qquad x\in A.
\end{equation}
  For the first sum in (\ref{G10}), we then obtain
\begin{eqnarray}
&&\sum_{w\in D(0,r-s)}~~\sum_{r< |y|\leq r+r^{1/(2+\bb)}}
G_A(x,y)p_1(y,w)\nn\\
&&\leq c\sum_{k=1}^{ r^{1/(2+\bb)}}~~\sum_{y\in U_k}
G_A(x,y)\frac{1}{(|y|-(r-s))^{3+2\beta}}\nn\\
&&
\leq c\sum_{k=1}^{ r^{1/(2+\bb)}} \frac{k}{(k-1+s)^{3+2\beta}}
\leq c\sum_{j=0}^{\ff} \frac{c}{(j+s)^{2+2\beta}}\nn\\
&&\leq \frac{c}{s^{1+2\beta}}.\label{G11}
\end{eqnarray}

For the second sum in (\ref{G10}), we use (\ref{p1.2j}) to bound $G_A(x,y)$ by
$c\log R$ and
  the fact that the cardinality of $U_k$ is bounded by $c(r+k)$ to obtain
\begin{eqnarray}
&&\sum_{w\in D(0,r-s)}~~\sum_{ r+r^{1/(2+\bb)}< |y|\leq R}
G_A(x,y)p_1(y,w)\nn\\
&&\leq c\sum_{k=r^{1/(2+\bb')}}^{\ff }~~\sum_{y\in U_k}
G_A(x,y)\frac{1}{(|y|-(r-s))^{3+2\beta}}\nn\\
&&
\leq c(\log R) \sum_{k=r^{1/(2+\bb)}}^{\ff } \frac{r+k}{k^{3+2\beta}}\nn\\
&&\leq c(\log R)
[r(r^{1/(2+\bb)})^{-(2+2\beta)}+(r^{1/(2+\bb)})^{-(1+2\beta)}].\label{G12}
\end{eqnarray}
By our assumptions on $R, r$ and $s$  we have our desired estimate, and it
only remains to prove (\ref{G6}).

We divide the proof of (\ref{G6}) into two steps.

\noindent{\bf Step 1.}  Suppose
$1\leq k\leq r^{1/(2+\bb)}$ and
$x\in D(0,r+k-2)$. Then there exists $0<c<1$ not depending
on
$k,r,$ or
$x$ such that
\begin{equation}\label{G1}
\P^x(T_{D(0,r)}< T_{D(0,r+k)})\geq c/k.
\end{equation}

\proof This is trivial if $k=1,2$ so we may assume that
$k\geq 3$ and that $x\in  D(0,r)^{ c}$.
By a rotation of the coordinate system, assume $x=(|x|,0)$. Let 
$Y_k=X_k\cdot (0,1)$, i.e., $Y_k$ is the second component of the 
random vector $X_k$.
Since the covariance matrix of $X$ is the identity, this is also true after
a rotation, so $Y_k$ is symmetric, mean 0, and the variance of $Y_1$ is 1.
Let $S_1$ be the line segment connecting $(r-1,-k^{1+\bb/2})$ with
$(r-1,k^{1+\bb/2})$ and $S_2$ the line segment connecting
$(r+k-1,-k^{1+\bb/2})$ to $( r+k-1,k^{1+\bb/2})$.
We have $S_1\subset D(0,r)$ because
$$(r-1)^2+(k^{1+\bb/2})^2=r^2-2r+1+k^{2+\bb}\leq r^2$$
by our assumption on $k$. Similarly
$S_2\subset D(0,r+k)$ because
$$(r+k-1)^2+(k^{1+\bb/2})^2=(r+k)^2-2(r+k)+1+k^{2+\bb}\leq (r+k)^2.$$
Let $L_i$ be the line containing $S_i$, $i=1,2$, and $Q$ the rectangle
whose left and right sides are $S_1$ and $S_2$, resp.

If $X$ hits $L_1$ before $L_2$ and does not exit $Q$ before hitting
$L_1$,  then
$T_{D(0,r)}<T_{D(0,r+k)}$. If $T_{L_1}\land T_{L_2}$ is less than
$k^{2+\bb/2}$ and $Y$ does not move more than $k^{1+\bb/2}$ is time
$k^{2+\bb/2}$, then
$X$ will not hit exit $Q$ through the top or bottom sides of $Q$.
Therefore
\begin{align}
\P^x(T_{D(0,r)}<T_{D(0,r+k)})&\geq \P^x(T_{L_1}<T_{L_2})-
\P^x(T_{L_1}\land T_{L_2}\geq k^{2+\bb/2})\nn\\
&~~~-\P^x(\max_{j\leq k^{2+\bb/2}} |Y_j|\geq k^{1+\bb/2})\nn\\
&=I_1-I_2-I_3. \label{G2}
\end{align}
By the gambler's ruin estimate \cite[(4.2)]{L3}, we have
\begin{equation}
\label{G3}
I_1\geq c\frac{(r+k-1)-|x|}{(r+k-1)-(r-1)}\geq \frac{c}{k}.
\end{equation}
It follows from the one dimensional case of (\ref{Lesc.10}) that for some
$\rho<1$ and all sufficiently large $k$
\begin{equation}
I_2\leq\P^x(T_{L_1}\land T_{L_2}\geq k^{2+\bb/2})\leq
\rho^{k^{\bb/4} }=o(1/k).\label{Lesc.10a}
\end{equation}
For $I_3$ we truncate the one-dimensional random walk $Y_j$ at level
$k^{1+\bb/4}$ and use Bernstein's inequality. If we let
$\xi_j=Y_j-Y_{j-1}$, $\xi'_j=\xi_j 1_{(|\xi_j|\leq k^{1+\bb/4})}$
and $Y'_j=\sum_{i\leq j} \xi'_i$, then
\begin{equation}
\label{G4}
\P^x(Y_j\ne Y'_j\mbox{ for some } j\leq k^{2+\bb/2})\leq \frac{k^{2+\bb}}
{(k^{1+\bb/4})^{3+2\beta}}=o(1/k).
\end{equation}
By Bernstein's inequality (\cite{Bennett})
\begin{equation}
\label{G5}
\P^x(\max_{j\leq k^{2+\bb/2}}|Y'_j|\geq k^{1+\bb/2})
\leq \exp\Big(-\frac{k^{2+\bb}}{2k^{2+\bb/2}+\frac23
k^{1+\bb/2}k^{1+\bb/4}}
\Big).
\end{equation}
This is also $o(1/k)$.

\noindent {\bf Step 2.} Let
$$J_k=\max_{x\in A} \sum_{y\in U_k} G_A(x,y).$$
By the strong Markov property, we see that the maximum is
taken when $x\in U_k$.
Also, by the strong Markov property,
\begin{equation}
\label{G8}
J_k\leq m+\sup_{x\in U_k}\E^x\Big[\sum_{y\in U_k} G_A(X_m,y)\Big].
\end{equation}
We have
\begin{align}
\E^x\Big[\sum_{y\in U_k} G_A(X_m,y);
&X_m\notin D(0,r+k-3)\Big]\nn\\
&\leq J_k\P^x(X_m\notin D(0,r+k-3))
\label{G1.9A}
\end{align}
and using (\ref{G1})
\begin{align}
\E^x\Big[\sum_{y\in U_k} G_A(X_m,y);
&X_m\in D(0,r+k-3)\Big]\nn\\
&\leq  \E^x[\E^{X_m}[J_k; T_{U_k}< T_{D(0,r)}];
X_m\in D(0,r+k-3)]\nn\\
&\leq J_k(1-c/k) \P^x(X_m\in D(0,r+k-3)).
\label{G1.9B}
\end{align}
Using (\ref{92.30}),
there exists $m$ and $\eps>0$ such that for all $x\in U_k$
\begin{equation}
\label{G7}
\P^x(X_m\in D(0,r+k-3))\geq \eps.
\end{equation}
Using  (\ref{G7}) and combining (\ref{G1.9A}) and (\ref{G1.9B})
  we see that  for all $x\in U_k$
$$\E^{ x}\Big[\sum_{y\in U_k} G_A(X_m,y)\Big]\leq J_k(1-c'/k)$$
for some $0<c'<1$.
Then by (\ref{G8})
\begin{equation}
\label{G9}
J_k\leq m+J_k(1-c'/k),
\end{equation}
and solving for $J_k$ yields (\ref{G6}).
\qed

Using (\ref{p1.2b}) and (\ref{G10a}) we can obtain that for some $\de>0$
\begin{eqnarray}
         &&
         \PPP^x\(T_{ D(0,r)}<T_{D(0,R)^c}\,;\,T_{ D(0,r)} =T_{\partial D(0,r)_{
s}}\)
\label{p1.2bu}\\
         &&
         = {\log (R/ |x|)\over
\log (R/ r)}+O( r^{-\de})+O(s^{ -1-2\beta}).  \nn
\end{eqnarray}

We now prove a  bound for the Green's function in the exterior of the disk
$D( 0,n)$.
\begin{lemma}\label{lem-gext}
\begin{equation} G_{ D( 0,n)^{ c}}( x,y)\leq c\log ( |x|\wedge |y|),\hspace{
.2in}x,y\in D( 0,n)^{ c}.\label{ext.1}
\end{equation}
\end{lemma}
\noindent {\bf Proof of Lemma~\ref{lem-gext}:} Since $G_{ D( 0,n)^{ c}}(
x,y)=\PPP^x\(T_{y}<T_{D(0,n)}\)G_{ D(  0,n)^{ c}}( y,y)$, and using the
symmetry of $G_{ D( 0,n)^{ c}}( x,y)$, it suffices to show that
\begin{equation} G_{ D( 0,n)^{ c}}( x,x)\leq c\log ( |x|),\hspace{ .2in}x\in D(
0,n)^{ c}.\label{ext.2}
\end{equation} Let
\begin{eqnarray} && U_1=0,\nn\\ && V_i=\min\{k>U_i: |X_k|< n
\mbox{ or } |X_k|\geq |x|^8\},\qquad  i=1,2,\ldots,\nn\\ &&
U_{i+1}=\min\{k>V_i: X_k=x\}, \qquad i=1,2,\ldots\nn
\end{eqnarray} Then using the strong Markov property and (\ref{p1.2j})
\begin{eqnarray} && G_{ D( 0,n)^{ c}}( x,x)=\E^x \lc\sum_{k<T_{ D( 0,n)}}
1_{\{x\}}(X_k)\rc\label{ext.3}\\ &&\leq
\sum_{i=1}^\infty
\E^x\Big[ \sum_{U_i\leq k< V_i} 1_{\{x\}}(X_k); U_i<T_{ D( 0,n)}\Big]\nn\\
&&\leq \sum_{i=1}^\infty \E^x[G_{_{ D( 0,|x|^8)}}(X_{U_i},x);  U_i<T_{ D(
0,n)}]\nn\\ &&\leq c_2 \log(|x|^8)
\sum_{i=1}^\infty \P^x(U_i<T_{ D( 0,n)}).\nn
\end{eqnarray} We have
$$\P^x(U_{i+1}<T_{ D( 0,n)})\leq \E^x[\P^{X_{U_i}}( T_{D( 0,|x|^8)^{ c}}<T_{
D( 0,n)}); U_i<T_{ D( 0,n)}].$$ By (\ref{p1.2bfa}),
\begin{eqnarray} &&
\P^{X_{U_i}}(T_{D( 0,|x|^8)^{ c}}<T_{ D( 0,n)})=\P^x(T_{D( 0,|x|^8)^{ c}}<T_{
D( 0,n)})\label{ext.4}\\ &&\leq
\frac{\log(|x|/n)+O(n^{-1/4})}{\log(|x|^8/n)}\nn\\ &&=\frac{\log |x|-\log
n+O(n^{-1/4})}{8\log |x|-\log n}\nn\\ &&\leq \frac{\log |x|+1}{7\log |x|}\leq
\frac27.\nn
\end{eqnarray} Therefore
$$\P^x (U_{i+1}<T_{ D( 0,n)})\leq \tfrac27 \P^x (U_i<T_{ D( 0,n)}).$$ By
induction $\P^x(U_i<T_{ D( 0,n)})\leq (\tfrac27)^i$, hence $\sum_i
\P^x(U_i<T_{ D( 0,n)})<c_4$. Together with (\ref{ext.3}), this proves
(\ref{ext.2}).
\qed

\begin{lemma}\label{lem-extjump}
\begin{equation}
\sup_{x\in D(0,n+s)^c}\PPP^{x}(T_{ D(0,n+s)}
\neq T_{\partial D(0,n)_{ s}})\leq cn^{ 2}\log (n)\,s^{-3-2\bb}+
cn^{-1-\bb}.\label{extjump.1}
\end{equation}
\end{lemma}

\noindent{\bf Proof of Lemma \ref{lem-extjump}:} Using (\ref{ext.1})
\begin{eqnarray}
         &&\sup_{x\in D(0,n+s)^c}\PPP^{x}(T_{ D(0,n+s)}
\neq T_{\partial D(0,n)_{ s}})
           \label{ext.60}\\
         &&=\sup_{x\in D(0,n+s)^c }
\sum_{\stackrel{y\in D(0,n+s)^c}{w\in  D(0,n)}} G_{D(0,n)^c }(x,y )p_{ 1}(y,w )
           \nn\\
         &&\leq c
\sum_{\stackrel{y\in D(0,n+s)^c}{w\in  D(0,n)}}\log( |y|)\, p_{ 1}(y,w )
\leq c\sum_{y\in D(0,n+s)^c}\log( |y|) \P(|X_{ 1}|\geq |y|-n )
           \nn\\
         &&
\leq c\sum_{y\in D(0,n+s)^c}\log( |y|)( |y|-n)^{ -3-2\bb}
           \nn\\
         &&
\leq c\log (n)\sum_{n+s\leq |y|\leq 2n }( |y|-n)^{ -3-2\bb} +c\sum_{2n<|y|
}\log( |y|)( |y|-n)^{ -3-2\bb}
           \nn\\
         &&
\leq cn^{ 2}\log (n)\,s^{-3-2\bb}+ cn^{-1-\bb}.
           \nn
\end{eqnarray}\qed

\section{Harnack inequalities}\label{sec-harn}

We next present some Harnack inequalities tailored to our needs.  We continue
to  assume that Condition A   and (\ref{1.0}) hold.

\begin{lemma}[Interior Harnack Inequality]\label{lem-H} Let $e^{ n}\leq r=
R/n^{ 3}$. Uniformly for $x,x'\in D(0,r)$ and
$y\in \partial D(0,R)_{ n^{ 4}}$
\begin{eqnarray}
         && H_{D(0,R)^c}(x,y)=\(1+O( n^{-3})\) H_{D(0,R)^c}(x',y).\label{p1.5}
\end{eqnarray}

Furthermore, uniformly in $x\in \partial D(0,r)_{ n^{ 4}}$ and
$y\in
\partial D(0,R)_{ n^{ 4}}$,
\bea &&
\PPP^{ x}(X_{ T_{D(0,R)^c }}=y\,,\,T_{D(0,R)^c }<T_{ D(0,r/n^{ 3}) } )
\label{p1.5h}\\ && =\(1+O( n^{-3})\)
\PPP^{ x}(T_{D(0,R)^c }<T_{ D(0,r/n^{ 3}) })H_{D(0,R)^c}(x,y).\nn
\eea
\end{lemma}

\noindent {\bf Proof of Lemma~\ref{lem-H}:}  It suffices to prove (\ref{p1.5})
with
$x'=0$. For any
$y\in \partial D(0,R)_{ n^{ 4}}$ we have the last exit decomposition
\begin{eqnarray}\qquad
          H_{D(0,R)^c}(x,y)&=&\sum_{z\in D(0,R)-D(0,3R/4)} G_{D(0,R)}(x,z)p_{
1}( z,y)\nn\\
         &+&\sum_{z\in D(0,3R/4)-D(0,R/2)} G_{D(0,R)}(x,z)p_{ 1}( z,y)\nn\\
         &+&
\sum_{z\in D(0,R/2)} G_{D(0,R)}(x,z)p_{ 1}( z,y).
         \label{p1.6}
\end{eqnarray} Let us first show  that uniformly in
$x\in D(0,r)$ and $z\in D(0,3R/4)-D(0,R/2)$
\begin{equation} G_{D(0,R)}(x,z)=\(1+O( n^{-3})\)G_{D(0,R)}(0,z).\label{b2.1}
\end{equation} To this end, note that by (\ref{p2.2}), uniformly in
$x\in D(0,r)$ and $y\in D(0,R/2)^{ c}$
\bea  a(y-x)&=&{ 2\over \pi}\log |x-y| +k +O( |x-y|^{ -1})\label{hq2.2}\\ &=&{
2\over \pi}\log |y| +k +O(n^{ -3})\nn\\ &=& a(y)+O(n^{ -3}).\nn
\eea Hence, by (\ref{p2.3}), and using the fact that $H_{D(0,R)^{ c} 
}(z,\cdot)$
is a probability, we see that uniformly in $x\in D(0,r)$ and
$z\in D(0,R)-D(0,R/2)$
\bea &&  G_{D(0,R)}(x,z) = G_{D(0,R)}(z,x)\label{hq2.3}\\ &&=\lc
\sum_{y\in D(0,R)^{ c}}H_{D(0,R)^{ c} }(z,y)a( y-x)\rc-a(z-x)\nn\\ &&=\lc
\sum_{y\in D(0,R)^{ c}}H_{D(0,R)^{ c}}(z,y) a( y)\rc- a( z)+O(n^{ -3})\nn\\
&&=G_{D(0,R)}(z,0)+O(n^{ -3}).\nn
\eea By (\ref{p1.2}), $G_{D(0,R)}(z,0)\geq c>0$ uniformly for $z\in
D(0,3R/4)-D(0,R/2)$, which completes the proof of (\ref{b2.1}).

We next show  that uniformly in
$x\in D(0,r)$ and $z\in D(0,R)-D(0,3R/4)$
\begin{equation} G_{D(0,R)}(x,z)=\(1+O(
n^{-3})\)G_{D(0,R)}(0,z).\label{b2.1b}
\end{equation} Thus, let
$z\in D(0,R)-D(0,3R/4)$, and use  the strong Markov property to see that
\begin{equation}G_{D(0,R)}(z,x) =
\E^{ z}\(  G_{D(0,R)}(X_{ T_{ D(0,3R/4)}},x)\,;\,T_{ D(0,3R/4)}<T_{ D(0,R)}
\)
\label{hq2.4}
\end{equation}
\[=\E^{ z}\(  G_{D(0,R)}(X_{ T_{ D(0,3R/4)}},x)\,;\,T_{ D(0,3R/4)}<T_{ D(0,R)}
\,;\,|X_{ T_{ D(0,3R/4)}}|>R/2\)\]
\[+\E^{ z}\(  G_{D(0,R)}(X_{ T_{ D(0,3R/4)}},x)\,;\,T_{ D(0,3R/4)}<T_{ D(0,R)}
\,;\,|X_{ T_{ D(0,3R/4)}}|\leq R/2\)\] By (\ref{p1.2j}) and 
(\ref{ext.60}) we can
bound the last term by
\begin{equation} c(\log R)\,\P^{ z}\(|X_{ T_{ D(0,3R/4)}}|\leq R/2\)\leq c(\log
R)^{ 2}R^{ -1-2\bb}.\label{hq2.47}
\end{equation} Thus we can write (\ref{hq2.4}) as
\bea &&  \E^{ z}\(  G_{D(0,R)}(X_{ T_{ D(0,3R/4)}},x)\,;\,T_{ D(0,3R/4)}<T_{
D(0,R)}
\,;\,|X_{ T_{ D(0,3R/4)}}|>R/2\)  \nn\\ &&\hspace{ 1in}=G_{D(0,R)}(z,x) +O(
R^{-1-\bb}).
\label{hq2.48}
\eea Applying (\ref{b2.1}) to the first line of (\ref{hq2.48}) and 
comparing the
result with (\ref{hq2.48}) for $x=0$ we see that uniformly in
$x\in D(0,r)$ and $z\in D(0,R)-D(0,3R/4)$
\begin{equation} G_{D(0,R)}(x,z)=\(1+O( n^{-3})\)G_{D(0,R)}(0,z)+O(
R^{-1-\bb}).\label{b2.1bx}
\end{equation} Since, by  (\ref{92.13}), $G_{D(0,R)}(0,z)\geq c/R$, we obtain
         (\ref{b2.1b}).

If $X$ has finite range then the last sum in (\ref{p1.6}) is zero and 
the proof of
(\ref{p1.5}) is complete. Otherwise, assume that (\ref{condA}) holds.  We next
show that
\begin{equation}
\sum_{z\in D(0,R/2)} G_{D(0,R)}(x,z)p_{ 1}( z,y)=O(R^{-3-\bb}).\label{p1.5j}
\end{equation} By (\ref{p1.2j}) we have $G_{D(0,R)}(x,z)=O( \log R)$ so that
the sum in (\ref{p1.5j})
         is bounded by $O( \log R)$ times the probability of a jump
         of size greater than $|y|-R/2\geq R/2$, which gives (\ref{p1.5j}).

To prove (\ref{p1.5}) it now suffices to show that uniformly in
$y\in \partial D(0,R)_{ n^{ 4}}$
\begin{equation} R^{-3-\bb}\leq cn^{ -3} H_{D(0,R)^c}(0,y).\label{b1.5ax}
\end{equation} Using once more the fact that $G_{D(0,R)}(0,z)\geq c/R$ by
(\ref{92.13}),
\bea
         H_{D(0,R)^c}(0,y)&=&\sum_{ z\in D( 0,R)}G_{D(0,R)}(0,z)p_{ 1}(
z,y)\label{b1.5bv}\\ &\geq & C \(\sum_{ z\in D(0,R)}p_{ 1}(y,z)\)R^{ -1}.\nn
\eea (\ref{b1.5ax}) then follows from (\ref{condA}) and our assumption that
$R\geq e^{ n}$,
         completing the proof of (\ref{p1.5}).

Turning to \req{p1.5h}, we have
\bea &&
\label{y1.5h}
\PPP^{ x}(X_{ T_{D(0,R)^c }}=y\,,\,T_{D(0,R)^c }<T_{ D(0,r/n^{ 3}) } )
\\ && =H_{ D(0,R)^c}(x,y)-
\PPP^{ x}(X_{ T_{D(0,R)^c }}=y\,,\,T_{D(0,R)^c }>T_{ D(0,r/n^{ 3}) } ).\nn
\eea By the strong Markov property at $T_{D(0,r/n^{ 3})}$
\bea &&
\PPP^{ x}(X_{ T_{D(0,R)^c }}=y\,,\,T_{D(0,R)^c }>T_{ D(0,r/n^{ 3}) }
)\label{y1.7h}\\ && =\E^{ x}(H_{D(0,R)^c}(X_{T_{ D(0,r/n^{ 3}) }},y)
;\,T_{D(0,R)^c  }>T_{ D(0,r/n^{ 3}) }).\nn
\eea By (\ref{p1.5}), uniformly in $w\in D(0,r/n^{ 3})$,
\[ H_{ D(0,n)^c}(w,y) =
\(1+O( n^{-3})\) H_{ D(0,n)^c}(x,y)\,.
\] Substituting back into (\ref{y1.7h}) we have
\beaa &&
\PPP^{ x}(X_{ T_{D(0,R)^c }}=y\,,\,T_{D(0,R)^c }>T_{ D(0,r/n^{ 3}) }
)\label{1.5tight}
\\ && =\(1+O( n^{-3})\)
\PPP^{ x}(T_{D(0,R)^c }>T_{ D(0,r/n^{ 3})})H_{D(0,R)^c}(x,y).\nn
\eeaa

Combining this with (\ref{y1.5h}) we obtain
\bea &&
\PPP^{ x}(X_{ T_{D(0,R)^c }}=y\,,\,T_{D(0,R)^c }<T_{ D(0,r/n^{ 3}) }
)\label{1.5tighter}
\\ && =\(\PPP^{ x}(T_{D(0,R)^c }<T_{ D(0,r/n^{ 3}) })+O( n^{-3})\)
H_{D(0,R)^c}(x,y).\nn
\eea Since, by (\ref{p1.2bfa})
\begin{equation}
\inf_{ x\in \partial D(0,r)_{  n^{ 4}}}
\PPP^{ x}(T_{D(0,R)^c }<T_{ D(0,r/n^{ 3}) })\geq 1/4,
\end{equation}
          we obtain \req{p1.5h} which completes the proof of Lemma~\ref{lem-H}.
\qed

Using the notation $n,r,R$ of the last Theorem, in preparation for 
the proof of
the exterior Harnack inequality we  now establish a uniform lower bound
      for the Greens function in the exterior of a disk:
\begin{equation} G_{ D( 0,r+n^{ 4})^{ c}}( x,y)\geq c>0 ,\hspace{
.2in}x\in\partial  D(0,R)_{ n^{4}},\,y\in D( 0,R)^{c}.\label{ext.7}
\end{equation}

Pick some $x_{ 1}\in \partial D(0,R)$ and proceeding clockwise choose points
$x_{ 1},\ldots, x_{ 36}\in \partial D(0,R)$ which divide
$\partial D(0,R)$ into
$36$ approximately equal arcs.  The distance between any two adjacent such
points is
$\sim 2R\sin (5^{ \circ})\approx .17R$. It then follows from (\ref{p1.2b}) that
for any
$j=1,\ldots, 36$
\bea &&
\inf_{ x\in T_{ D(x_{ j},R/5)}}\P^{ x}( T_{ D(x_{ j+1},R/5)}<T_{ D(0,r+n^{
4})})\label{ext.20}\\ &&\geq \inf_{ x\in T_{ D(x_{ j+1},2R/5)}}\P^{ x}( T_{
D(x_{ j+1},R/5)}<T_{ D(x_{ j+1},R/2)^{ c}})\geq c_{ 1}>0\nn
\eea for some $c_{ 1}>0$ independent of $n,r,R$ for $n$ large and where we
set $x_{ 37}=x_{ 1}$. Hence, a simple argument using the strong Markov
property shows that
\begin{equation}
\inf_{ j,k}\inf_{ x\in T_{ D(x_{ j},R/5)}}\P^{ x}( T_{ D(x_{k},R/5)}<T_{
D(0,r+n^{ 4})})\geq c_{ 2}=:c^{ 36}_{ 1}\label{ext.21}
\end{equation} Furthermore, it follows from (\ref{p1.2c}) that for any
$j=1,\ldots, 36$
\bea &&
\inf_{ x,x'\in T_{ D(x_{ j},R/5)}}\P^{ x}( T_{ x'}<T_{ D(0,r+n^{
4})})\label{ext.22}\\ &&\geq \inf_{ x\in T_{   D(x',2R/5)}}\P^{ x}( T_{ x'}<T_{
D(x',R/2)^{ c}})\geq c_{ 3}/\log R\nn
\eea for some $c_{ 3}>0$ independent of $n,r,R$ for $n$ large. Since
$\partial  D(0,R)_{R/100}\subseteq \cup_{ j=1}^{ 36} D(x_{ j},R/5)$,
combining (\ref{ext.21}) and (\ref{ext.22})  we see that
\begin{equation}
\inf_{ x,x'\in \partial  D(0,R)_{R/100}}\P^{ x}( T_{ x'}<T_{ 
D(0,r+n^{ 4})})\geq
c_{ 4}/\log R.\label{ext.23}
\end{equation} It then follows from (\ref{p2.5}) that
\begin{eqnarray} &&
\inf_{ x,x'\in \partial  D(0,R)_{ R/100}}G_{ D( 0,r+n^{ 4})^{
c}}(x,x')\label{ext.24}\\ && =\inf_{ x,x'\in \partial D(0,R)_{R/100}}\P^{ x}(
T_{ x'}<T_{ D(0,r+n^{ 4})})G_{ D( 0,r+n^{ 4})^{ c}}(x',x')
\nonumber\\ &&
\geq (c_{ 4}/\log R)\, G_{ D( x',R/2)}(x',x')\geq c_{ 5}>0
\nonumber
\end{eqnarray} for some $c_{ 5}>0$ independent of $n,r,R$ for
$n$ large.

Using the strong Markov property, (\ref{ext.24}), and (\ref{ext.60}) 
we see that
\bea &&\inf_{ z\in D( 0,1.01R)^{ c},\,x\in \partial  D(0,R)_{ R/100}} G_{ D(
0,r+n^{ 4})^{ c}}(z,x)\label{ext.25}\\ &&
\geq \E^z[G_{ D( 0,r+n^{ 4})^{ c}}(X_{T_{D( 0,1.01R)}},x); X_{T_{D(
0,1.01R)}}\in
\partial  D(0,R)_{R/100}]\geq c>0.\nn
\eea Together with (\ref{ext.24}) this  completes the proof of (\ref{ext.7}).

We next observe that uniformly in $x\in\partial  D(0,R)_{ n^{ 4}}\,,\,z\in
D(0,2r)-D(0,5r/4)$,  by (\ref{ext.7}) and (\ref{p1.2bfa}),
\begin{eqnarray} \qquad&& G_{D(0,r+n^{ 4})^{ c}}(x,z)=G_{D(0,r+n^{ 4})^{
c}}(z,x)\label{ext.35}\\ && =\E^{ z}\lc  G_{D(0,r+n^{ 4})^{ c}}(X_{T_{D( 0,R)^{
c}} }, x)\,;\,  T_{D( 0,R)^{ c}}<T_{D(0,r+n^{ 4})}
\rc
\nonumber\\ &&\geq c \P^{ z}\lc  T_{D( 0,R)^{c}}<T_{D(0,r+n^{ 4})}\rc
\nonumber\\ &&\geq c /\log n.
\nonumber
\end{eqnarray}

Our final observation is  that for any $\ep>0$, uniformly in
$x\in\partial  D(0,R)_{ n^{ 4}}\,,\,z\in D(0,2r)-D(0,r+( 1+\ep)n^{ 4})$,
\begin{equation} G_{D(0,r+n^{ 4})^{ c}}(x,z)\geq c(R\log R)^{ 
-1}.\label{ext.50}
\end{equation} To see this, we first use the strong Markov property together
with (\ref{ext.23}) to see that uniformly in
$x, x'\in
\partial  D(0,R)_{ n^{ 4}}$ and
$z\in D(0,2r)-D(0,r+( 1+\ep)n^{ 4})$,
\bea G_{D(0,r+n^{ 4})^{ c}}(x,z)&\geq &\P^{ x}( T_{ x'}<T_{ D(0,r+n^{
4})})G_{D(0,r+n^{ 4})^{ c}}(x',z)\label{ext.52}\\ &\geq &cG_{D(0,r+n^{ 4})^{
c}}(x',z) /\log R.\nn
\eea In view of (\ref{92.13}),  if $x'\in
\partial  D(0,R)_{ n^{ 4}}$ is chosen as close as possible to the 
ray from the the
origin which passes through $z$
\begin{equation} G_{D(0,r+n^{ 4})^{ c}}(x',z)\geq G_{D(x',|x'|-(r+n^{
4}))}(x',z)\geq c R^{ -1}.\label{ext.51}
\end{equation} Combining the last two displays proves (\ref{ext.50}).

\begin{lemma}[Exterior Harnack Inequality]\label{lem-HT} Let $e^{ n}\leq r=
R/n^{ 3}$. Uniformly for $x,x'\in
\partial  D(0,R)_{ n^{ 4}}$ and
$y\in \partial  D(0,r)_{ n^{ 4}}$
\begin{equation} H_{D(0,r+n^{ 4})}(x,y)=
\(1+O(n^{ -3}\log n)\)H_{D(0,r+n^{ 4})} (x',y).\label{p1.5t}
\end{equation}

   Furthermore,  uniformly in $x\in \partial D(0,R)_{ n^{ 4}}$ and
$y\in  \partial  D(0,r)_{ n^{ 4}}$,
\begin{eqnarray} &&\qquad
\PPP^{ x}(X_{ T_{D(0,r+n^{ 4})}}=y\,;\,T_{ D(0,r+n^{ 4})}<T_{D(0,n^{ 3}R)^{
c}})
\label{p1.5br}\\ && =\Big(1+O(n^{ -3}\log n)\Big)
\PPP^{x}(T_{ D(0,r+n^{ 4})}<T_{D(0,n^{ 3}R)^{ c}} )H_{D(0,r+n^{ 4})}(x,y),\nn
\end{eqnarray}
         and uniformly in $x,x'\in  \partial D(0,R)_{ n^{ 4}}$ and $y\in
\partial   D(0,r)_{ n^{ 4}}$,
\begin{eqnarray} &&\qquad
\PPP^{ x}(X_{ T_{ D(0,r+n^{ 4})}}=y\,;\,T_{ D(0,r+n^{ 4})}<T_{D(0,n^{ 3}R)^{
c}} )
\label{p1.5b}\\ && =\Big(1+O(n^{ -3}\log n)\Big)
\PPP^{ x'}(X_{ T_{ D(0,r+n^{ 4})}}=y\,;\,T_{ D(0,r+n^{ 4})}<T_{D(0,n^{ 3}R)^{
c}} ).\nn
\end{eqnarray}
\end{lemma}

{\bf Proof of Lemma~\ref{lem-HT}:}  For any
$x\in
\partial  D(0,R)_{ n^{ 4}}$ and
$y\in \partial  D(0,r)_{ n^{ 4}}$ we have the last exit decomposition
\begin{eqnarray}\qquad
          H_{D(0,r+n^{ 4})}(x,y)&=&\sum_{z\in D(0,5r/4)-D(0,r+n^{ 4})}
G_{D(0,r+n^{ 4})^{ c}}(x,z)p_{ 1}( z,y)\nn\\
         &+&
\sum_{z\in D(0,2r)-D(0,5r/4)} G_{D(0,r+n^{ 4})^{ c}}(x,z)p_{ 1}( z,y)\nn\\
         &+&
\sum_{z\in D^{ c}(0,2r)} G_{D(0,r+n^{ 4})^{ c}}(x,z)p_{ 1}( z,y).
         \label{q1.6}
\end{eqnarray}

Let us first show  that uniformly in
$x,x'\in
\partial  D(0,R)_{ n^{ 4}}$ and $z\in D(0,2r)-D(0,5r/4)$
\begin{equation} G_{D(0,r+n^{ 4})^{ c}}(x,z)=\(1+O( n^{-3}\log
n)\)G_{D(0,r+n^{ 4})^{ c}}(x',z).\label{q2.1}
\end{equation} To this end, note that by (\ref{p2.2}), uniformly in
$x,x'\in
\partial  D(0,R)_{ n^{ 4}}$ and $y\in D(0,2r)$
\bea  a(x-y)&=&{ 2\over \pi}\log |x-y| +k +O( |x-y|^{ -1})\label{q2.2}\\ &=&{
2\over \pi}\log R +k +O(n^{ -3})\nn\\ &=& a(x'-y)+O(n^{ -3}),\nn
\eea and with $N\geq n^{ 3}R$ the same result applies with $y\in D(0,N)^{ c}$.
Hence, by (\ref{p2.3}) applied to the finite set
$A(r+n^{ 4}, N )=:D(0,N)-D(0,r+n^{ 4})$, and using the fact that
$H_{A(r+n^{ 4}, N )^{ c} }(z,\cdot)$ is  a probability, we see that 
uniformly in
$x,x'\in
\partial  D(0,R)_{ n^{ 4}}$ and $z\in D(0,r+n^{ 4})^{ c}$
\bea && \qquad G_{A(r+n^{ 4}, N )}(x,z) = G_{A(r+n^{ 4}, N
)}(z,x)\label{q2.3}\\ &&=\lc \sum_{y\in D(0,r+n^{ 4})\cup D(0,N)^{
c}}H_{A(r+n^{ 4}, N )^{  c} }(z,y)a( y-x)\rc-a(z-x)\nn\\ &&=\lc
\sum_{y\in D(0,r+n^{ 4})\cup D(0,N)^{  c}}H_{A(r+n^{ 4}, N )^{ c}}(z,y) a(
y-x')\rc\nn\\ && \hspace{ 2in}- a( z-x')+O(n^{ -3})\nn\\ &&=G_{A(r+n^{ 4}, N
)}(z,x')+O(n^{ -3}).\nn
\eea Since this is uniform in $N\geq n^{ 3}R$, using (\ref{ext.1}) we can apply
the dominated convergence theorem  as
$N\rar
\ff$ to see that  uniformly in
$x,x'\in
\partial  D(0,R)_{ n^{ 4}}$ and $z\in D(0,2r)$
\begin{equation} G_{D(0,r+n^{ 4})^{ c}}(x,z)=G_{D(0,r+n^{ 4})^{
c}}(z,x')+O(n^{  -3}).\label{ext.42}
\end{equation} Applying (\ref{ext.35}) now establishes (\ref{q2.1}).

We show next that uniformly in
$x,x'\in
\partial  D(0,R)_{ n^{ 4}}$ and $z\in D(0,5r/4)-D(0,r+n^{ 4})$
\begin{equation} G_{D(0,r+n^{ 4})^{ c}}(x,z)=\(1+O( n^{-3}\log
n)\)G_{D(0,r+n^{ 4})^{  c}}(x',z)+O(r^{ -1-2\bb}\log r).\label{ext.43}
\end{equation} To see this we use the strong Markov property together with
(\ref{q2.1})  and (\ref{ext.1}) to see that
\begin{eqnarray} && G_{D(0,r+n^{ 4})^{ c}}(x,z)=G_{D(0,r+n^{ 4})^{
c}}(z,x)\label{ext.44}\\ && =\E^{ z}\lc  G_{D(0,r+n^{ 4})^{ c}}(X_{D(0,5r/4)^{
c}  },x)\,;\,T_{D(0,5r/4)^{ c} } <T_{D(0,r+n^{ 4}) }\rc,\nonumber
\end{eqnarray} and this in turn is bounded by
\bea &\E^{ z}\lc  G_{D(0,r+n^{ 4})^{ c}}(X_{D(0,5r/4)^{ c}
},x)\,;\,T_{D(0,5r/4)^{ c} } <T_{D(0,r+n^{ 4}) }\,,\,|X_{D(0,5r/4)^{ c} }|\leq
2r\rc\nn\\ &+\E^{ z}\lc  G_{D(0,r+n^{ 4})^{ c}}(X_{D(0,5r/4)^{ c}
},x)\,;\,T_{D(0,5r/4)^{ c} } <T_{D(0,r+n^{ 2}) }\,,\,|X_{D(0,5r/4)^{ c} }|>
2r\rc\nn\\ &\leq \(1+O( n^{-3}\log n)\)
\E^{ z}\lc  G_{D(0,r+n^{ 2})^{ c}}(X_{D(0,5r/4)^{ c} 
},x')\,;\,T_{D(0,5r/4)^{ c}
}  <T_{D(0,r+n^{ 4}) }\rc\nn\\ &\hspace{.3in}+c\log (R)\, \P^{ z}\lc
|X_{D(0,5r/4)^{ c} }|> 2r\rc.\label{S31}\nn
\eea
   Using (\ref{ext.44}) we see that
$$\E^{ z}\lc  G_{D(0,r+n^{ 2})^{ c}}(X_{D(0,5r/4)^{ c} },x')\,;\,T_{D(0,5r/4)^{
c} }  <T_{D(0,r+n^{ 4}) }\rc =G_{D(0,r+n^{ 4})^{ c}}(x',z).$$
         As in (\ref{jumpo.2}), by using (\ref{p1.2j}) and (\ref{2.0}) we have
\begin{eqnarray} &&
\P^{ z}\lc  |X_{D(0,5r/4)^{ c} }|> 2r\rc   \label{jumpo.2a}\\
         &&=
\sum_{\stackrel{|y|<5r/4}{2r<|w|}} G_{D(0,5r/4) }(z,y )p_{ 1}(y,w )
           \nn\\
         &&\leq c \log r
\sum_{|y|<5r/4}
\P(|X_{ 1}|\geq 3r/4 )
           \nn\\
         &&\leq c \log r
\sum_{|y|<5r/4} { 1\over |r |^{ 3+2\bb}}
\leq cr^{-1-2\bb}\log r.
           \nn
\end{eqnarray} This establishes (\ref{ext.43}).

         We next show that
\begin{equation}
\sum_{z\in D^{ c}(0,2r)} G_{D(0,r+n^{ 4})^{ c}}(x,z)p_{ 1}(
z,y)=O(r^{-3-\bb}).\label{q1.5j}
\end{equation}
         It follows from (\ref{ext.1}) that
$G_{D(0,r+n)^{ c}}(x,z)= O( \log R)$, so that the sum in (\ref{p1.5j})
         is bounded by $O( \log R)$ times the probability of a jump
         of size greater than $|z|-r-n\geq r-n\geq r/2$, which gives 
(\ref{q1.5j}).

To prove (\ref{p1.5t}) it thus suffices to show that uniformly in
$x\in
\partial  D(0,R)_{ n^{ 4}}$ and
$y\in \partial  D(0,r)_{ n^{ 4}}$
\begin{equation} r^{ -1-2\bb}\log r\leq cn^{ -3}  H_{D(0,r+n^{
4})}(x,y).\label{b1.5a}
\end{equation} Using  the last exit decomposition together with (\ref{ext.50})
we obtain
\bea &&
          H_{D(0,r+n^{ 4})}(x,y)=\sum_{z\in D(0,r+n^{ 4})^{ c}} G_{D(0,r+n^{
4})^{ c}}(x,z)p_{ 1}( z,y)\label{b1.5s}\\ &&\geq c(R\log R)^{ -1}
\sum_{r+(1+\ep)n^{ 4}\leq |z|\leq 2r} p_{ 1}( z,y)\nn
\eea for any $\ep>0$. Note that the annulus $\{ z\,|\, r+(1+\ep)n^{ 4}\leq
|z|\leq 2r\}$ contains the disc $D(v, 2(1+\ep)n^{ 4} )$ where
$v=(r+3(1+\ep)n^{ 4} )y/|y|$, and we have $2(1+\ep)n^{ 4}\leq |y-v|\leq
3(1+\ep)n^{ 4}$.  Thus
\begin{eqnarray} \sum_{r+(1+\ep)n^{ 4}\leq |z|\leq 2r} p_{ 1}( z,y)&\geq &
\sum_{z\in D(v, 2(1+\ep)n^{ 4} )} p_{ 1}( z,y)\label{b1.5sx}\\&= &
\sum_{z\in D(v, 2(1+\ep)n^{ 4} )} p_{ 1}( z-v,y-v)
\nonumber\\&= &
\sum_{z\in D(0, 2(1+\ep)n^{ 4} )} p_{ 1}( z,y-v).
\nonumber
\end{eqnarray} Hence by (\ref{condA}) and our assumption that
$r\geq e^{ n}$
\begin{equation}
          H_{D(0,r+n^{4})}(x,y)\geq c(R\log R)^{ -1}e^{ -( 1+\ep)^{ 1/4}\bb n}
\geq cr^{ - 1-\ep-( 1+\ep)^{ 1/4}\bb},\label{ext.64}
\end{equation} and thus (\ref{b1.5a}), and hence (\ref{p1.5t}), follows by
taking
$\ep$ small.

The rest of Lemma~\ref{lem-HT} follows as in the proof of
Lemma~\ref{lem-H}.\qed

\section{Local time estimates and upper bounds}\label{sec-hit}
\label{sec-upperbound}

          The following simple lemma, the analog of \cite[Lemma 
2.1]{jhung}, will
be used repeatedly.

\begin{lemma}\label{lem-hit}  For $|x_0|<R$,
\begin{equation}
\E^{x_0}(L^{ 0}_{T_{D(0, R)^{ c}}})=   G_{D(0,R)}(x_0,0).\label{tohelet}
\end{equation} For all
$z\geq 1$
\begin{equation}
    \P^{x_0}( L^{ 0}_{T_{D(0, R)^{ c}}}\geq zG_{D(0,R)}(0,0))\leq c\sqrt{ z}e^{
-z}\,\label{prob}
\end{equation} for some $c<\ff$ independent of $x_0, z, R$.

Let $x_0\neq 0$.   Let $0<\varphi\leq 1$ and set
$\la=\varphi/G_{D(0,R)}(0,0)$. Then
\begin{eqnarray} &&
\E^{x_0}\(e^{ -\la L^{ 0}_{T_{D(0, R)^{ c}}}}\)\label{babe.1}\\ &&=1-{\log
({R\over |x_0|})+O( |x_{ 0}|^{-1/4})
\over
\log (R)}{\varphi\over 1+\varphi }\(1+O({ 1\over \log (|x_0|)} )\)\nn
\end{eqnarray}
\end{lemma}

\noindent {\bf Proof of Lemma~\ref{lem-hit}:}  Since
\begin{equation} \label{c1.1}  L^{ 0}_{T_{D(0, R)^c}}=\sum_{i<T_{D(0,R)^c}}
1_{\{X_i=0\} },
\end{equation}  (\ref{tohelet}) follows from (\ref{p1.1}). Then we have by the
strong Markov property that
\begin{eqnarray}
\E^{x_0}(L^{ 0}_{T_{D(0, R)^{ c}}})^k & \!\!= \!\!& k!
\E^{x_0}\(\sum_{0 \leq j_1 \leq \cdots\leq j_k \leq T_{D(0, R)^{ c}}}
\prod_{i=1}^k 1_{\{X_{j_{ i} }=0\} } \) \nn \\ &\!\!= \!\!& k!
\E^{x_0}\(\sum_{0 \leq j_1 \leq \cdots\leq j_{ k-1} \leq
         T_{D(0, R)^{ c}}}
\prod_{i=1}^{ k-1} 1_{\{X_{j_{ i} }=0\} } G_{D(0,R)}(0,0)\)
\nn \\ &\!\!=\!\!& k\E^{x_0}(L^{ 0}_{ T_{ D(0,R )^c}})^{ k-1}G_{D(0,R)}(0,0).
\,,
\nn
\end{eqnarray} By induction on
$k$,
\begin{equation}
\label{secmom}
\E^{x_0}(L^{ 0}_{T_{D(0, R)^{ c}}})^k= k!G_{D(0,R)}(x_0,0)(
G_{D(0,R)}(0,0))^{k -1}.
\end{equation}
     To prove (\ref{prob}), use (\ref{secmom}),  (\ref{mark.1})  and Chebyshev
to obtain
\begin{equation}
\label{prob1} \P^{x_0}( L^{ 0}_{T_{D(0, R)^{ c}}}\geq zG_{D(0,R)}(0,0))\leq {
k!\over z^{ k}}\,
\end{equation} then take $k=[z]$ and use Stirling's formula.

For (\ref{babe.1}), note that, conditional  on  hitting $0$,  $L^{ 0}_{T_{D(0,
R)^{ c}}}$ is a geometric random variable with  mean
$G_{D(0,R) }(0,0)$.   Hence,
\bea &&
    \E^{x_0}\(e^{ -\la  L^{ 0}_{T_{D(0, R)^{ c}}}}\)
\nn
\\ && =1-\P^{ x_0}\(T_{ 0}<T_{D(0, R)^{c}}\)\nn\\ &&\hspace{ .5in}
    +\P^{ x_0}\(T_{ 0}<T_{D(0, R)^{c}}\)\({ 1\over
    (e^{ \la }-1 )G_{D(0,R) }(0,0)+1 }
\).\label{8.15}
\eea

Since by (\ref{p2.5})
\begin{equation} { 1\over G_{D(0,R) }(0,0)}=O(1/
\log (R))\label{8.15a}
\end{equation} we have
\begin{equation} (e^{ \la }-1 )G_{D(0,R) }(0,0)+1 =1+\varphi +O\Big({ 1\over
\log (R)} \Big)\label{8.15aa}
\end{equation} and (\ref{babe.1}) then follows from (\ref{8.15}) and
(\ref{p1.2c}).
\qed

We next provide the required upper bounds in Theorem
\ref{theo-late}. Namely,  we will show that for any $a \in (0,2]$
\begin{equation}
\label{3.1b}
\limsup_{m\rar \ff}{ \log \Big|\Big\{x\in  D(0, m):\; L^{ x}_{  T_{D(0, m)^{
c}}}\geq (2a/\pi) (\log m)^2\Big\}\Big|\over \log m}\leq 2-a
\hspace{.2in}a.s.
\end{equation}

To see this fix $\ga>0$ and note that by (\ref{prob})  and 
(\ref{p2.5}), for some
$0<\de<\ga$, all $x\in D(0, m)$ and all large enough $m$
\begin{equation}
\PPP^{0}\(  \frac{L^{ x}_{ T_{D(x, 2m)^{ c}}}}{(\log m)^2}\geq 2a/\pi
\)\leq m^{ -a+\de}\label{prob2}
\end{equation} Therefore
\begin{eqnarray} &&
\PPP^{0}\(\Big|\Big\{x\in  D(0, m):\;
         \frac{L^{ x}_{ T_{D(0, m)^{ c}}}}{(\log m)^2}\geq 2a/\pi
\Big\}\Big|\geq m^{2-a+\ga}\)\label{3.1ba}\\ &&\leq
m^{-(2-a)-\ga}\E^{0}\(\Big|
\Big\{x\in  D(0, m):\;
         \frac{L^{ x}_{ T_{D(0, m)^{ c}}}}{(\log m)^2}\geq 2a/\pi
\Big\}\Big|\)\nn\\ &&=m^{-(2-a)-\ga}\sum_{x\in  D(0, m)}\PPP^{0}\(
\frac{L^{ x}_{ T_{D(0, m)^{ c}}}}{(\log m)^2}\geq 2a/\pi
\)\nn\\  &&\leq m^{-(2-a)-\ga}\sum_{x\in  D(0, m)}\PPP^{0}\(
\frac{L^{ x}_{ T_{D(x, 2m)^{ c}}}}{(\log m)^2}\geq 2a/\pi \)\nn\\
         &&\leq m^{-(\ga-\de )}.\nn
\end{eqnarray} Now apply our result to $m=m_{ n}=e^{ n}$ to see by
Borel-Cantelli that for some $N(\om)<\ff$ a.s. we have that for all
$n\geq N(\om)$
\begin{equation}
         \Big| \Big\{x\in  D(0, e^{ n}):\;\frac{L^{ x}_{ T_{D(0, e^{ 
n})^{ c}}}}{(\log
e^{ n})^2}\geq 2a/\pi \Big\}\Big|\leq e^{(2-a+\ga)n}.\label{prob3}
\end{equation} Then if $e^{ n}\leq m\leq e^{ n+1}$
\begin{eqnarray} &&
         \Big| \Big\{x\in  D(0, m):\;\frac{L^{ x}_{ T_{D(0, m)^{ c}}}}{(\log
m)^2}\geq 2a/\pi
\Big\}\Big|\label{prob4}\\ &&
         \leq \Big| \Big\{x\in  D(0, e^{ n+1}):\;\frac{L^{ x}_{ 
T_{D(0, e^{ n+1})^{
c}}}}{(\log e^{ n})^2}\geq 2a/\pi
\Big\}\Big|\nn\\ && = \Big| \Big\{x\in  D(0, e^{ n+1}):\;\frac{L^{ 
x}_{ T_{D(0,
e^{ n+1})^{ c}}}}{(\log e^{ n+1})^2}\geq 2a( 1+1/n)^{- 2}/\pi
\Big\}\Big|\nn\\ &&
         \leq e^{(2-a( 1+1/n)^{ -2}+\ga)(n+1)}\leq m^{(2-a( 1+1/n)^{
-2}+\ga)(n+1)/n}.\nn
\end{eqnarray}
        (\ref{3.1b}) now follows by first letting $n\to \infty$ and then letting
        $\ga\rar 0$.
\qed

\section{Lower bounds for probabilities}\label{sec-lowerbound}

Fixing $a<2$, we prove in this section
\begin{equation}
\label{p.1}
\liminf_{m\rar \ff}{ \log \Big|\Big\{x\in  D(0, m):\; L^{ x}_{  T_{D(0, m)^{
c}}}\geq (2a/\pi) (\log m)^2\Big\}\Big|\over \log m}\geq 2-a
\hspace{.2in}a.s.
\end{equation} In view of (\ref{3.1b}), we will obtain Theorem
\ref{theo-late} .

We start by constructing a subset of the set appearing in  (\ref{p.1}), the
probability of which is easier to bound below. To this end fix $n$, and let
$r_{ n,k}=e^{ n} n^{ 3(n-k)},\,k=0,\ldots, n$. In particular,  $r_{ 
n,n}=e^{ n}$
and
$r_{ n,0}=e^{ n} n^{ 3n}$. Set
$K_n=16 r_{ n,0}= 16 e^{ n}n^{ 3 n}$.

Let $U_{ n}=[2r_{ n,0},3r_{ n,0}]^{ 2}\subseteq D(0, K_n)$. For
$x\in U_{ n}$, consider the $x$-bands
$\partial D(x,r_{n,k})_{ n^{ 4}};\,k=0,\ldots,n$. We use the abbreviation
$r'_{ n,k}=r_{ n,k}+n^{ 4}$. For $x\in U_{n}$ we will say that the path {\bf
does not skip
$x$-bands} if
\begin{itemize}
\item[(1)] $T_{D(x,r'_{n,0})}<T_{ D(0,K_n)^{ c}}$ and
$T_{D(x,r_{n,0}')}=T_{\partial D(x,r_{n,0})_{n^{4}} }$.
\item[(2)] For any
$t<T_{ D(0,K_n)^{ c}}$ such that $X_{ t}\in \partial D(x,r_{n,k})_{ n^{ 4}}$ we
have:
\item[(2$'$)] if $k=0$ then
\[\(T_{D(x,r_{n,1}')}\wedge T_{D(0,K_n)^{ c}}\)\circ\theta_{ t}=
\(T_{\partial D(x,r_{n,1})_{ n^{ 4}}}
\wedge T_{D(0,K_n)^{ c}}\)\circ\theta_{ t},\]
\item[(2$''$)] if $k=1,\ldots,n-1$ then
\[\(T_{D(x,r_{n,k+1}')}\wedge T_{D(x,r_{n,k-1})^{ c}}\)\circ\theta_{ t}=
\(T_{\partial D(x,r_{n,k+1})_{ n^{ 4}}}
\wedge T_{\partial D(x,r_{n,k-1})_{ n^{ 4}}}\)\circ\theta_{ t},\]
\item[(2$'''$)] if $k=n$ then
\[\(T_{D(x,r_{n,n-1})^{ c}}\)\circ\theta_{ t}=
\(T_{\partial D(x,r_{n,n-1})_{ n^{ 4}}}\)\circ\theta_{ t}.\]
\end{itemize}

        For $x
\in D(0,K_n)$, let $N_{n,k}^x$ denote the number of excursions from
$D(x,r_{n,k-1})^{ c}$ to
$D(x,r_{n,k}')$ until time $T_{D(0,K_n)^{ c}}$. Set
$\NNN_k=3ak^2\log k$, and $k_{ 0}=4\vee\inf\{ k\,|\,\NNN_k\geq 2k\}$. We
will say that a point
$x\in U_{ n}$ is {\bf $n$-successful} if the path does not skip
$x$-bands, $N^x_{n,k}=1\,,\,\forall k=1,\ldots,k_{ 0}-1$ and
\begin{equation} \NNN_k-k\leq N^x_{n,k}\leq
\NNN_k+k\hspace{.3in}\forall k=k_{ 0},\ldots,n.
\label{pmperf}
\end{equation}

         Let
$\{Y(n,x)\,;\,x\in U_{ n}\}$ be the collection of random variables defined by
\[Y(n,x)=1\hspace{.1in}\mbox{if $x$ is $n$-successful}\] and
$Y(n,x)=0$ otherwise. Set
$\bar{q}_{n,x}=\PPP(Y(n,x)=1)=\E(Y(n,x))$.

    The next lemma relates the notion of
$n$-successful and local times. As usual we write $\log_{ 2} n$ for
$\log\,\log n$.
\begin{lemma}\label{plem1} Let
\[\mathcal{S}_n=\{x\in U_{ n}\,|\,\mbox{$x$ is
$n$-successful}\}.\] Then for some $N( \om)<\ff$ a.s., for all $n\geq  N(
\om)$ and all $x\in \mathcal{S}_n$
\[
\frac{L^{ x}_{ T_{D(0, K_n)^{ c}}}}{ (\log K_n)^2}\geq  2a/\pi-2/\log_{ 2} n.\]
\end{lemma}

\noindent {\bf Proof of Lemma~\ref{plem1}:} Recall that if $x$ is n-successful
then
$N^x_{n,n}\geq \NNN_{n}-n=a(3n^2\log n)-n$.  Let $L^{ x,j}$ denote the
number of visits to $x$ during the
$j^{th}$ excursion from $\partial D(x,r_{n,n})_{ n^{ 4}}$ to $ D(x, 
r_{n,n-1})^{
c}$. Then
    for any $0<\la<\ff$
\begin{eqnarray} P_x &:=&\PPP\(L^{ x}_{ T_{D(0, K_n)^{ c}}}
\leq (2a/\pi-2/\log_{ 2} n)(\log K_n)^2
\,, \,x\in \mathcal{S}_n\)
\nn\\ &\leq &\PPP\(\sum_{j=0}^{\NNN_{n}-n}L^{ x,j}
\leq (2a/\pi-1/\log_{ 2} n)(3n\log n)^2
\)\nn\\ & \leq & \exp\(\la (2a/\pi-1/\log_{ 2} n)(3n\log n)^2\)
    E\(e^{ -\la \sum_{j=0}^{\NNN_{n}-n}L^{ x,j}} \).\label{pc3.6}
\end{eqnarray}
    If $\tau$ denotes the first time that the $(\NNN_{n}-n)^{th}$ excursion from
$  D(x, r_{n,n-1})^{ c}$ reaches
$\partial D(x,r_{n,n})_{ n^{ 4}}$ then by the strong Markov property
\begin{eqnarray} &&\qquad \E\(e^{ -\la \sum_{j=0}^{\NNN_{n}-n}L^{ x,j}}
\)\label{pc3.6a}\\ &&=\E\(e^{ -\la
\sum_{j=0}^{\NNN_{n}-n-1}L^{ x,j}}\E^{X_{\tau }}\(e^{ -\la L^{ x}_{ T_{D(x,
r_{n,n-1})^{ c}}}}\) \).\nn
\end{eqnarray} Set $\la=\phi/G_{D(x, r_{n,n-1}) }(x ,x)$. By (\ref{babe.1}),
with $r=r_{n,n}=e^{ n}, R=r_{n,n-1}=n^{3 }e^{ n}$, for any
$0<\phi\leq 1$ and large $n$
\begin{equation}
\sup_{y\in  \partial D(x,r_{n,n})_{ n^{ 4}}}\E^{y}\(e^{ -\la L^{ x}_{ T_{D(x,
r_{n,n-1})^{ c}}}}\)\leq \exp\(-{\(1-1/2\log n\)\phi \over 
1+\phi}\,3( \log n)/n
\).\label{pc3.6b}
\end{equation} Hence by induction
\begin{equation} \E\(e^{ -\la \sum_{j=0}^{\NNN_{n}-n} L^{ x,j} } \)
\leq \exp\(-{( 1-1/\log n)\phi \over 1+\phi}9an( \log n)^{ 2}
\).\label{pc3.6c}
\end{equation} Then with this choice of $\la$, noting that $G_{D(x, r_{n,n-1})
}(x ,x)\sim { 2\over
\pi}n$ by (\ref{p2.5}), we have
\begin{equation} P_x\leq \inf_{\phi>0 }\exp\(\Big\{\phi( 1-1/2\log_{ 2} n) -
{( 1-1/\log n)\phi \over 1+\phi}\Big\}9an( \log n)^{ 2}\).\label{pc3.6d}
\end{equation} A straightforward computation shows that
\begin{equation}
\inf_{\phi>0 }\(\phi\al -{ \phi\over 1+\phi
}\bb\)=-\(\sqrt{\bb}-\sqrt{\al}\)^{ 2}\label{8.21}
\end{equation} which is achieved for $\phi =\sqrt{\bb}/\sqrt{\al}-1$. Using
this in (\ref{pc3.6d}) we find that
\begin{equation} P_x\leq \exp\(-cn(\log n/\log_{ 2} n)^{ 2}\).\label{pc3.6e}
\end{equation} Note that $|U_{ n}|\leq e^{ cn\log n}$.
    Summing over all $x\in U_{ n}$ and then over $n$ and applying
Borel-Cantelli will then complete the proof of Lemma~\ref{plem1}.
\qed

The next lemma, which provides estimates for the first and second moments of
$Y(n,x)$, will be proved in the following sections. Recall $\bar q_{n,x}=
        \PPP(x\,
\mbox{ is $n$-successful})$.  Let $Q_n=\inf_{x\in U_n} \bar q_{n,x}$.
\begin{lemma}\label{pmomlb}  There exists $\de_n \to 0$ such that for all
$n \geq 1$,
\begin{equation}\qquad Q_{ n}
        \geq K_{n}^{-(a+\de_n)},
\label{pmomlb.1}
\end{equation} and
\begin{equation} Q_n\geq c \bar{q}_{n,x}
\label{pmomlb.1j}
\end{equation}
         for some $c>0$ and all $n$ and $x\in U_{ n}$.

There exists $C<\ff$ and $\de'_n \to 0$ such that for all $n$, $x
\neq y$ and
$l(x,y)=\min \{m\,:
\,D(x,r_{n,m})\cap D(y,r_{n,m})=\emptyset\} \leq n$
\begin{equation}
\E(Y(n,x)Y(n,y))\leq C Q_n^2  ( l(x,y)!)^{3a+\de'_{l(x,y) }}\;.
\label{pm3.5}
\end{equation}
\end{lemma}

\medskip\noindent{\bf Proof of Theorem \ref{theo-late}.} In view of
(\ref{3.1b}) we need only consider the lower bound. To prove (\ref{p.1}) we
will show that for any
$\de>0$ we can find $p_{ 0}>0$ and $N_{ 0}<\ff$ such that
\begin{equation}
\PPP^{0}\(\sum_{x\in U_{ n}}1_{\{ Y(n,x)=1\}}\geq K_{ n}^{2-a-\de}\)
\geq p_{0}\label{pm3.6}
\end{equation} for all $n\geq N_{ 0}$. Lemma \ref{plem1} will then imply that
for some $p_{ 1}>0$ and $N_{ 1}<\ff$
\begin{equation} \label{pm3.6z}
\PPP^{0}\(\Big|\Big\{x\in  D(0, K_{ n}):\; L^{ x}_{ T_{D(0, K_n)^{ c}}}\geq
(2a/\pi-2/\log n)  (\log K_{ n})^2\Big\}\Big|\geq K_{ n}^{2-a-\de}\)
\geq p_1.
\end{equation} for all $n\geq N_{ 1}$. As in the  proof of (\ref{3.1b}) and
readjusting $\de>0$ we can find $p_{ 2}>0$ and
$N_{ 2}<\ff$ such that
\begin{equation}
\PPP^{0}\(\Big|\Big\{x\in  D(0, n):\; L^{ x}_{ T_{D(0, n)^{ c}}}\geq (2a/\pi )
(\log n)^2\Big\}\Big|\geq n^{2-a-\de}\)\geq p_{2}\label{pm3.6zy}
\end{equation} for all $n\geq N_{ 2}$. Then by Lemma
\ref{LRatio}, with a further readjustment of $\de>0$ we will have that
\begin{equation}
\PPP^{0}\(\Big|\Big\{x\in  \Z^{ 2}:\; L^{ x}_{ n}\geq (a/2\pi ) (\log
n)^2\Big\}\Big|\geq n^{1-a/2-\de}\)\geq p_{3}\label{pm3.6ze}
\end{equation} for some  $p_{ 3}>0$ and all $n\geq N_{ 3}$ with
$N_{ 3}<\ff$. This estimate leads to (\ref{p.1}) as in the proof of 
Theorem 5.1 of
\cite{DPRZ4}.

Recall the Paley-Zygmund inequality (see
\cite[page 8]{Kahane}): for any $W\in L^{ 2}(
\Om)$ and $0<\la<1$
\begin{equation}
         \PPP(W\geq \la \E( W) )\geq ( 1-\la)^{ 2}{ (\E( W))^{ 2}\over
\E( W^{ 2})}.\label{PZ}
\end{equation} We will apply this with $W=W_{  n}=\sum_{x\in U_{ n}}1_{\{
Y( n,x)=1\}}$.  We see by (\ref{pmomlb.1}) of Lemma
\ref{pmomlb} that for some sequence
$\de_n\rar 0$
\begin{equation}
\E\(\sum_{x\in U_{ n}}1_{\{ Y(n,x)=1\}}\)=\sum_{x\in U_{ n}}
\bar{q}_{n,x}\geq K_n^{2-a-\de_n}.\label{pm3.7}
\end{equation}  Thus to complete the proof of (\ref{pm3.6}) it suffices to show
\begin{equation}
\E\(\lc\sum_{x\in U_{ n}}1_{\{ Y(n,x)=1\}}\rc^2\)\leq c\lc
\E\(\sum_{x\in U_{ n}}1_{\{ Y(n,x)=1\}}\)\rc^2\label{p.2}
\end{equation} for some $c<\ff$ all $n$ sufficiently large. Furthermore, using
(\ref{pm3.7}) it suffices to show that
\begin{equation}
\E\(\sum_{\stackrel{x,y\in U_{ n}}{x\neq y}}1_{\{ Y(n,x)=1\}}1_{\{
Y(n,y)=1\}}\)\leq c\lc \E\(\sum_{x\in U_{ n}}1_{\{
Y(n,x)=1\}}\)\rc^2\label{p.3}
\end{equation}

We let $C_m$ denote generic finite constants that are  independent of
$n$. The definition of $l(x,y)
\geq 1$ implies that
$|x-y|\leq 2r_{n,l(x,y)-1}$. Recall that because $n\geq l$, there are at most
$C_0 r_{ n,l-1}^2\leq C_0K^{ 2}_n/n^{ 6( l-1)}\leq C_0K^{ 2}_n\,l^{ 6}( l!)^{
-6}$
         points
$y$ in the ball of radius $2r_{n, l-1}$ centered at $x$. Thus it  follows from
Lemma
\ref{pmomlb} that
\begin{eqnarray} &&
\sum_{\stackrel{x,y\in U_{ n}}{2r_{ n,n}\leq |x-y|\leq 2r_{
n,0}}}\E\(Y(n,x)Y(n,y)\)\label{p.5a}\\ &&
          \leq C_1 \sum_{\stackrel{x,y\in U_{ n}}{2r_{ n,n}\leq 
|x-y|\leq 2r_{ n,0}}}
Q_n^2 ( l(x,y)!)^{3a+\de'_{l(x,y) }}\nn\\  &&\leq C_2 Q_n^2
\sum_{x\in U_n} \sum_{j=1}^n \sum_{\{y\, |\,  l(x,y)=j\}}
(j!)^{3a+\delta'_j}\nn\\ && \leq C_3 Q_n^2 K^{ 2}_n
\sum_{j=1}^{n} K^{ 2}_n\,j^{ 6}( j!)^{ -6}( j!)^{3a+\de'_{j }}\nn\\  &&
\leq C_4 (K^{ 2}_n Q_n)^2  \sum_{j=1}^{n}  j^{ 6}(j!)^{-3(2-a)+\de'_{j }}\nn\\
         &&\leq C_5 (K^{ 2}_n Q_n)^2\leq C_6\lc \E\(\sum_{x\in U_{
n}}Y(n,x)\)\rc^2
\nn
\end{eqnarray} where we used the fact from (\ref{pmomlb.1}) that
\begin{equation} K^{ 2}_n Q_n\leq c \sum_{x\in U_{ n}}\bar{q}_{n,x}
=c\E\(\sum_{x\in U_{ n}}Y(n,x)\).\label{p.5f}
\end{equation} Because $Y(n,x)\leq 1$ and $\E Y(n,x)=\bar q_{n,y}\leq cQ_n$,
we have
\begin{eqnarray} &&
\sum_{\stackrel{x,y\in U_{ n}}{|x-y|\leq 2r_{ n,n}}}\E\(Y(n,x)Y(n,y)\)
\leq \sum_{\stackrel{x,y\in U_{ n}}{|x-y|\leq 2r_{n,
n}}}\E\(Y(n,y)\)\label{p.5}\\&&
          \leq C_7 \sum_{x\,;\,|x|\leq 2r_{n, n}} K^{ 2}_nQ_n
          \leq C_8 \,e^{ 2n} K^{ 2}_nQ_n.
\nn
\end{eqnarray} By (\ref{pmomlb.1})
\[ K_n^2 Q_n\geq K_n^{2-a-\delta_n}\geq ce^{2n}.\] This and (\ref{p.5f})
show that the right hand side of  (\ref{p.5}) is bounded by
\begin{equation}\label{5.18}
        C_9\lc \E\(\sum_{x\in U_{ n}}Y(n,x)\)\rc^2.
\end{equation} We know that if $x,y\in U_n$, then $|x-y|\leq 2r_{n,0}$.
         Thus combining (\ref{p.5a}), (\ref{p.5}), and (\ref{5.18}) 
completes the
proof of (\ref{p.3}) and hence of (\ref{p.1}).
\qed

\medskip

\noindent{\bf Proof of Theorem \ref{theo-et}.}    The lower bound is an
immediate consequence of Theorem
\ref{theo-late}. The upper bound is a consequence of (\ref{prob}) as follows;
cf.~the proof of
\cite[(2.8)]{DPRZ4}.  Let $\delta>0$. By (\ref{p2.5}) and (\ref{prob})
\begin{eqnarray} && \P^0(\sup_{x\in D(0,R)} L^x_{T_{D(0,R)^c}}>
\frac{4}{\pi}(1+\delta) \log^2 R)
\label{pfub}\\ &\leq & \sum_{x\in D(0,R)}
\P^0(L^x_{T_{D(x,2R)^c}}>
\frac4{\pi}(1+\delta)
\log^2 R)\nn\\ &\leq & cR^2 ((1+2\delta)\log R)^{1/2} e^{-2(1+\delta/2)
\log R}
\leq c R^{-\delta/4}\nn
\end{eqnarray} for $R$ large. By Borel-Cantelli there exists
$M_0(\omega)$ such that if $m\geq M_0$, then
$$L^*_{T_{D(0,2^m)^c}}\leq \frac4{\pi} (1+\delta) \log^2 (2^m).$$ If
$m\geq M_0$, $2^m\leq n\leq 2^{m+1}$, and $m$ is large,
$$L^*_{T_{D(0,n)^c}}\leq L^*_{T_{D(0,2^{m+1})^c}}
\leq \frac{4}{\pi} (1+\delta) \log^2 (2^{m+1}) \leq \frac4{\pi} (1+2\delta)
\log^2 n.$$ Since $\delta$ is arbitrary, this and Lemma \ref{LRatio} prove the
upper bound.
\qed

\section{First moment estimates}\label{51momest}

\noindent {\bf Proof of  (\ref{pmomlb.1}) and (\ref{pmomlb.1j})}:  For
$x\in U_{ n}$ we begin by getting bounds on the probability that
$T_{D(x,r'_{n,0})}<T_{  D(0,K_n)^{ c}}$ and
$T_{D(x,r_{n,0}')}=T_{\partial D(x,r_{n,0})_{ n^{ 4}} }$. Since
\begin{eqnarray} &&
         \PPP \(T_{D(x,r'_{n,0})}<T_{ D(0,K_n)^{
c}}\,;\,T_{D(x,r_{n,0}')}=T_{\partial D(x,r_{n,0})_{ n^{ 4}} }\)\label{pr.1}
\\ &&
\geq  \PPP \(T_{D(x,r'_{n,0})}<T_{D(x,{1 \over 2}K_{ n})^{ c}}\,;
\,T_{D(x,r_{n,0}')}=T_{\partial D(x,r_{n,0})_{ n^{ 4}} }
\)\nn
\end{eqnarray}  we see from (\ref{p1.2bu})  that uniformly in $n$ and
$x\in U_{ n}$
\begin{equation}
\PPP\(T_{D(x,r'_{n,0})}<T_{ D(0,K_n)^{ c}}\,;\,T_{D(x,r_{n,0}')}=T_{\partial
D(x,r_{n,0})_{ n^{ 4}} }\) \geq c\label{pr.2}
\end{equation}  for some $c>0$. And since for $x\in U_{ n}$ and
$y\in
\partial D(x,r_{n,0})_{ n^{ 4}}$
\begin{eqnarray} &&
\PPP^{ y} \(T_{D(x,r'_{n,1})}<T_{ D(x,{1 \over 2}K_n)^{
c}}\,;\,T_{D(x,r_{n,1}')}=T_{\partial D(x,r_{n,1})_{ n^{ 4}} }\)\label{pr.3}\\
&&\hspace{ .3in}\leq
         \PPP^{ y} \(T_{D(x,r'_{n,1})}<T_{ D(0,K_n)^{
c}}\,;\,T_{D(x,r_{n,1}')}=T_{\partial D(x,r_{n,1})_{ n^{ 4}} }\)\nn\\ &&\leq
\PPP^{ y} \(T_{D(x,r'_{n,1})}<T_{ D(x,2K_n)^{
c}}\,;\,T_{D(x,r_{n,1}')}=T_{\partial D(x,r_{n,1})_{ n^{ 4}} }\)
\nn
\end{eqnarray}  we see from (\ref{p1.2bu}) that uniformly in $n$,
$x\in U_{ n}$ and $y\in
\partial D(x,r_{n,0})_{ n^{ 4}}$
\begin{equation}
         c/\log n\leq  \PPP^{ y} \(T_{D(x,r'_{n,1})}<T_{ D(0,K_n)^{ c}}\,;
\,T_{D(x,r_{n,1}')}=T_{\partial D(x,r_{n,1})_{ n^{ 4}} }\) \leq c'/\log
n.\label{pr.4}
\end{equation} Similarly, since for $x\in U_{ n}$ and $y\in
\partial D(x,r_{n,0})_{ n^{ 4}}$
\begin{eqnarray}
         \PPP^{ y} \(T_{D(0,K_{ n})^{ c}}<T_{  D(x,r'_{n, 1})}\)
\geq   \PPP^{ y} \(T_{ D(x,2K_{ n})^{ c}}<T_{D(x,r'_{n, 1})}\)
\label{pr.3a}
\end{eqnarray} we see from (\ref{p1.2bd}) that uniformly in $n$,
$x\in U_{ n}$ and $y\in
\partial D(x,r_{n,0})_{ n^{ 4}}$
\begin{equation}
         \PPP^{ y} \(T_{D(0,K_{ n})^{ c}}<T_{  D(x,r'_{n, 1})}\) \geq c>0.
\label{pr.5}
\end{equation}

These bounds will be used for excursions at the `top' levels. To  bound
excursions at `intermediate' levels we note that  using (\ref{p1.2bd}), we have
uniformly for
$x\in
\partial D(0,r_{n,l})_{ n^{ 4}}$, with
$1\leq l\leq n-1$
\begin{eqnarray}
         &&
         \PPP^x\(T_{D(0,r_{n, l-1})^{ c}}<T_{  D(0,r'_{ n,l+1})}\,;
\,T_{D(0,r_{n, l-1})^{ c}}=T_{\partial D(0,r_{n, l-1})_{ n^{ 4}}}\)
\label{ps.3a}\\ &&
         =1/2+O(n^{ -4-4\bb}),\nn
\end{eqnarray} and   using  (\ref{p1.2bu}), we have uniformly for
$x\in
\partial D(0,r_{n,l})_{ n^{ 4}}$, with
$1\leq l\leq n-1$
\begin{eqnarray}
         &&
         \PPP^x\(T_{D(0,r'_{n, l+1})}<T_{ D(0,r_{ n,l-1})^{ c}}\,;
\,T_{D(0,r'_{n, l+1})}=T_{\partial D(0,r_{n, l+1})_{ n^{ 4}}}\)
\label{ps.3b}\\ &&
         =1/2+O(n^{ -4-4\bb}).\nn
\end{eqnarray}

For excursions at the `bottom' level, let us  note, using an analysis 
similar to
that of  (\ref{c2.11k}), that uniformly in
$z\in D(0,r_{n,n})_{ n^{ 4}}$
\begin{equation}
\PPP^z\(T_{D(0,r_{n,n-1})^c} =T_{\partial D(0,r_{n, n-1})_{ n^{ 4}}}\)
=1+O(n^{ -4-4\bb} ).\label{ps.3c}
\end{equation}

Let $\bar{m}=(m_{ 2}, m_{ 3},\ldots, m_{ n})$
         and set $|\bar{m}|=2\sum_{ j= 2}^{ n}m_{ j}+1$. Let
$\mathcal{H}_{ n}(\bar{m})$,
         be the collection of maps, (`histories'),
\[\varphi: \{0,1,\ldots,  |\bar{m}|  \}\mapsto \{0,1,\ldots,n \}\] such that
$\varphi(0)=1,\,\varphi( j+1)=\varphi( j)\pm 1,\,|\bar{m}|
=\inf\{j\,;\,\varphi(j)=0
\}$
         and the number of upcrossings from $\ell-1$ to
$\ell$
\[u( \ell)=:|\{( j,j+1) \,|\, ( \varphi(j),\varphi(  j+1))=(\ell-1,\ell)
\}|=m_{\ell }.\]

Note that we cannot have any upcrossings from
$\ell$ to $\ell+1$ until we have first had an upcrossing from
$\ell-1$ to $\ell$. Hence the number of ways to partition the $u(
\ell+1)$ upcrossings from
$\ell$ to
$\ell+1$ among and after the $u( \ell)$  upcrossings from $\ell-1$ to
$\ell$ is  the same as the number of ways to partition
$u(\ell+1)$ indistinguishable objects into $u( \ell)$ parts, which is
\begin{equation} { u( \ell+1)+u( \ell)-1\choose u(
\ell)-1}.\label{44.1}
\end{equation}
         Since $u( \ell)=m_\ell$ and the mapping $\varphi$ is completely
determined once
         we know the relative order of all its upcrossings
\begin{equation} |\mathcal{H}_{ n}(\bar{m})|=
\prod_{\ell= 2}^{ n-1}{m_{ \ell+1}+m_{ \ell}-1\choose m_{
\ell}-1}.\label{44.2}
\end{equation}

        What we do next is estimate the probabilities of all possible 
orderings of
visits to
$\{\partial D(x,r_{n,j})_{n^4}: j=0, \ldots, n\}$. Let $\Om_{ x,n}$ denote the
set of random walk paths which  do not  skip $x$-bands
         until completion of the first
         excursion from
$ \partial D(x,r_{n,1})_{ n^{ 4}}$ to $
\partial D(x,r_{n,0})_{ n^{ 4}}$.
         To each random walk path $\om\in \Om_{ x,n}$ we assign a `history'
$h(\om )$ as follows.
         Let $\tau( 0)$ be the time of the first visit to
$\partial D(x,r_{n, 1})_{ n^{ 4}}$, and  define $\tau( 1), \tau( 2),
\ldots$ to be the successive hitting times of different elements of
\[\{\partial D(x,r_{n,0})_{ n^{ 4}},
\ldots,
\partial D(x,r_{n, n})_{ n^{ 4}}\}\] until the first downcrossing from
$\partial D(x,r_{ n,1})_{ n^{ 4}}$ to
$\partial D(x,r_{ n,0})_{ n^{ 4}}$. Setting $\Phi( y)=k$ if $y\in
\partial D(x,r_{ n,k})_{ n^{ 4}}$,
         let
$h(\om )( j)=\Phi( \om (\tau( j) ))$.  Let $h_{|_{k} }$ be the restriction of
$h$ to
$\{0,\ldots,k  \}$.  We claim that uniformly for any
         $\varphi\in \mathcal{H}_{ n}(\bar{m})$ and $z\in \partial 
D(x,r_{n,1})_{
n^{ 4}}$
\begin{eqnarray}
         &&
          \P^{ z}\lc h_{|_{ |\bar{m}|} }=\varphi\,;\,\Om_{ x,n}\rc
         =\( { 1\over 2} \)^{|\bar{m}|-m_{ n} }\lc 1+O(n^{-4-4\bb})
         \rc^{|\bar{m}| }.\label{44.3}
\end{eqnarray} To see this, simply use the strong Markov property successively
at the times
\[\tau( 0), \tau( 1),
\ldots, \tau(|\bar{m}|-1)\] and then use (\ref{ps.3a})-(\ref{ps.3c}).

Writing $m\stackrel{k}{\sim} \NNN_{k}$ if $m=1$ for $k<k_{ 0}$ and
$|m-\NNN_k| \leq k$ for $k\geq k_{ 0}$ we see
         that uniformly in $m_{ n}\stackrel{ n}{\sim}\NNN_{ n}$ we have that
$\lc 1+O(n^{-4-4\bb})
         \rc^{|\bar{m}| }=1+O(n^{-1-3\bb})$. Combining this with 
(\ref{44.2}) and
(\ref{44.3}) we see that  uniformly in $z\in
\partial D(x,r_{n,1})_{ n^{ 4}}$
\begin{eqnarray} &&
\sum_{\stackrel{m_{ 2}, \ldots, m_{ n}} {m_{\ell}\stackrel{\ell}{\sim}\NNN_{
\ell}} }
         \P^{ z}\lc h_{|_{ |\bar{m}|} }\in
\mathcal{H}_{ n} (\bar{m})\,;\,\Om_{ x,n}\rc\label{44.5m}\\ &&
         =( 1+O(n^{-1-3\bb}))\sum_{\stackrel{m_{ 2}, \ldots, m_{ n}}
{m_{\ell}\stackrel{\ell}{\sim}\NNN_{
\ell}} }\( { 1\over 2} \)^{|\bar{m}|-m_{ n} }
\prod_{\ell= 2}^{ n-1}{m_{ \ell+1}+m_{ \ell}-1\choose m_{
\ell}-1} \nonumber\\ &&
         =( 1+O(n^{-1-3\bb})){ 1\over 4}\sum_{\stackrel{m_{ 2}, \ldots, m_{ n}}
{m_{\ell}\stackrel{\ell}{\sim}\NNN_{
\ell}} }
\prod_{\ell= 2}^{ n-1}{m_{ \ell+1}+m_{ \ell}-1\choose m_{
\ell}-1}\( { 1\over 2} \)^{m_{ \ell+1}+m_{ \ell} }. \nonumber
\end{eqnarray} Here we used the fact that $m_{ 2}=1$ so that
$|\bar{m}|-m_{ n} =2+\sum_{\ell= 2}^{ n-1}(m_{ \ell+1}+m_{
\ell})$.

\begin{lemma}\label{lem-stirling} For some $C=C(a) <\infty$ and all
$k\geq 2$,
$|m-\NNN_{k+1}| \leq k+1$, $|\ell+1-\NNN_{k}| \leq  k$,
\begin{equation} { C^{ -1}k^{-3a-1}\over \sqrt{\log k}}
\leq {m+\ell\choose \ell}\( { 1\over 2} \)^{m+\ell+1}
\leq { Ck^{-3a-1}\over \sqrt{\log k}}.\label{44.5}
\end{equation}
\end{lemma}

\noindent {\bf Proof of Lemma \ref{lem-stirling}}: It suffices to consider
$k \gg 1$ in which case the binomial coefficient in (\ref{44.5}) is well
approximated by Stirling's formula
\[ m!=\sqrt{2\pi}m^m e^{-m}\sqrt{m} (1+o(1)) \;.
\] With $\NNN_k = 3a k^2 \log k$ it follows that for some
$C_1<\ff$ and all
$k$ large enough, if $|m-\NNN_{k+1}| \leq 2 k$,
$|\ell-\NNN_{k}| \leq 2 k$ then
\begin{equation} |\frac{m}{\ell} - 1 - \frac{2}{k}| \leq \frac{C_1}{k
\log k} \;.
\label{m1.3la}
\end{equation} Hereafter, we use the notation
         $f \sim g$ if $f/g$ is bounded and bounded away from zero as
$k \to
\infty$, uniformly in
$\{m:\,|m-\NNN_{k+1}| \leq 2 k\}$ and $\{\ell:|\ell-\NNN_{k}|
\leq 2 k\}$. We then have by the preceding observations that
\begin{equation} {m+\ell\choose \ell}\( { 1\over 2} \)^{m+\ell+1}
\sim
\frac{(m+\ell)^{m+\ell}}{\sqrt{\ell}\, \ell^\ell m^m}
\( { 1\over 2} \)^{m+\ell} \sim \frac{\exp(-\ell I( {m \over
\ell}))}{\sqrt{k^2 \log k}} \;,
\label{Lkest}
\end{equation} where
$$ I(\la)= - (1+\la) \log (1+\la) + \la \log \la + \la \log 2 + \log 2 \;.
$$ The function $I(\la)$ and its first order derivative vanishes at
$1$, with the second derivative $I_{\la \la} (1) = 1/2$. Thus, by a Taylor
expansion to second order of $I(\la)$ at $1$, the estimate 
(\ref{m1.3la}) results
with
\begin{equation} | I({m \over \ell}) - \frac{1}{k^2} | \leq
\frac{C_2}{k^2
\log k}
\label{Ipest}
\end{equation} for some $C_2 < \infty$, all $k$ large enough and
$m,\ell$ in the range considered here.  Since
$|\ell - 3a k^2 \log k| \leq 2k$, combining (\ref{Lkest}) and (\ref{Ipest}) we
establish (\ref{44.5}).
\qed

Using the last Lemma we have that
\begin{eqnarray} &&
\sum_{\stackrel{m_{ 2}, \ldots, m_{ n}} {m_{\ell}\stackrel{\ell}{\sim}\NNN_{
\ell}} }
\prod_{\ell= 2}^{ n-1}{ C^{ -1}\ell^{-3a-1}\over \sqrt{\log
\ell}}\nn\\ &&
\hspace{ .5in}\leq
\sum_{\stackrel{m_{ 2}, \ldots, m_{ n}} {m_{\ell}\stackrel{\ell}{\sim}\NNN_{
\ell}} }
\prod_{\ell= 2}^{ n-1}{m_{ \ell+1}+m_{ \ell}-1\choose m_{
\ell}-1}\( { 1\over 2} \)^{m_{ \ell+1}+m_{ \ell} }\label{44.5n}\\ &&
\hspace{ 1in}\leq
\sum_{\stackrel{m_{ 2}, \ldots, m_{ n}} {m_{\ell}\stackrel{\ell}{\sim}\NNN_{
\ell}} }
\prod_{\ell= 2}^{ n-1}{ C\ell^{-3a-1}\over \sqrt{\log \ell}}.
\nonumber
\end{eqnarray} Using the fact that $|\{
m_{\ell}\,|\,m_{\ell}\stackrel{\ell}{\sim}\NNN_{
\ell}\}|=2\ell+1$, this shows that for some $C_1<\ff$,
\begin{eqnarray} && n\prod_{\ell= 2}^{ n-1}{ C_1^{ -1}\ell^{-3a}\over
\sqrt{\log \ell}}\nn\\ &&
\hspace{ .5in}\leq
\sum_{\stackrel{m_{ 2}, \ldots, m_{ n}} {m_{\ell}\stackrel{\ell}{\sim}\NNN_{
\ell}} }
\prod_{\ell= 2}^{ n-1}{m_{ \ell+1}+m_{ \ell}-1\choose m_{
\ell}-1}\( { 1\over 2} \)^{m_{ \ell+1}+m_{ \ell} }\label{44.5p}\\ &&
\hspace{ 1in}\leq n\prod_{\ell= 2}^{ n-1}{ C_1\ell^{-3a}\over
\sqrt{\log \ell}}.  \nonumber
\end{eqnarray}

Since for any $c<\ff$, for some $\zeta_{ n},\zeta'_{ n}\rar 0$
\begin{equation} nc^{ n}\prod_{\ell= 2}^{ n-1}\log
\ell=n^{n\zeta_{ n} }=( n!)^{\zeta'_{ n} }
\label{44.5q}
\end{equation} we see that  for some $\de_{1, n}, \de_{ 2,n}\rar 0$
\begin{equation}
\sum_{\stackrel{m_{ 2}, \ldots, m_{ n}} {m_{\ell}\stackrel{\ell}{\sim}\NNN_{
\ell}} }
\prod_{\ell= 2}^{ n-1}{m_{ \ell+1}+m_{ \ell}-1\choose m_{
\ell}-1}\( { 1\over 2} \)^{m_{ \ell+1}+m_{ \ell} }
         =( n!)^{-3a-\de_{1, n}}=r_{n,0}^{-a-\de_{2, n}}.\label{44.6q}
\end{equation}

(\ref{pr.2})-(\ref{pr.5}) and (\ref{44.5m})   show that for some
$0<c,c'<\ff$
\begin{eqnarray} && {c \over \log n}\sum_{\stackrel{m_{ 2},
\ldots, m_{ n}} {m_{\ell}\stackrel{\ell}{\sim}\NNN_{
\ell}} }
\prod_{\ell= 2}^{ n-1}{m_{ \ell+1}+m_{ \ell}-1\choose m_{
\ell}-1}\( { 1\over 2} \)^{m_{ \ell+1}+m_{ \ell} }\label{44.6f}\\ &&
\hspace{ .3in}\leq Q_{ n}=\inf_{ x\in U_{ n}} \PPP(x\,
\mbox{ is $n$-successful})\leq \sup_{ x\in U_{ n}} \PPP(x\,
\mbox{ is $n$-successful})\nonumber\\ &&
\hspace{ .6in}\leq {c' \over \log n}
\sum_{\stackrel{m_{ 2}, \ldots, m_{ n}} {m_{\ell}\stackrel{\ell}{\sim}\NNN_{
\ell}} }
\prod_{\ell= 2}^{ n-1}{m_{ \ell+1}+m_{ \ell}-1\choose m_{
\ell}-1}\( { 1\over 2} \)^{m_{ \ell+1}+m_{ \ell} }.\nonumber
\end{eqnarray}

Together with (\ref{44.6q}) this gives (\ref{pmomlb.1}) and
(\ref{pmomlb.1j}).\qed

In the remainder of this section we prove two lemmas needed to complete the
proof of Lemma \ref{pmomlb}.

Let $\Om_{ x,n,i,m}^{i-1,\ldots,j}$ denote the set of random walk paths which
do not skip
$x$-bands on excursions between levels $k=i-1,i,\ldots,j$
         until completion of the first
$m$ excursions from
$D(x,r'_{n,i})$ to $ D(x,r_{n,i-1})^{ c}$ and let $N_{n,i,m,k}^x$ denote the
number of excursions from
$  D(x,r_{n,k-1})^{ c}$ to $
         D(x,r'_{n,k})$  until completion of the first
$m$ excursions from
$ D(x,r'_{n,i})$ to $
         D(x,r_{n,i-1})^{ c}$.

\begin{lemma}\label{cond-corl} We can find
$C<\ff$ and $\de_{3, l}\rar 0$ such that for all $n$ and $1\leq l< n$,
\bea &&
         \sum_{\stackrel{m_{ l}, \ldots, m_{ n}} {m_{k}\stackrel{k}{\sim}\NNN_{
k}} }\PPP
\( N_{n,k}^x =m_{k},\; k=l+1,\ldots,n\,;\,\Om_{ 
x,n,l,m_{l}}^{l-1,\ldots,n}\,|\,
N_{n,l}^x =m_{l}\) \label{44.9qd}\\ &&\hspace{ 1in}
\leq C\, Q_n  \,\,( l!)^{3a+\de_{3,l }}.\nn
\eea
\end{lemma}

\noindent{\bf Proof of Lemma \ref{cond-corl}:} The analysis of this section
shows that uniformly in
$m_{k}\stackrel{k}{\sim}\NNN_{k},\; k=l,l+1,\ldots,n$ and
$z\in \partial D(x,r_{n,l})_{ n^{ 4}}$
\begin{eqnarray}
         &&
         \PPP^{ z} \( N_{n,l,m_{ l},k}^x =m_{k},\; k=l+1,\ldots,n\,;\,\Om_{
x,n,l,m_{ l}}^{l-1,\ldots,n}\)
\label{44.8}\\
         &&
         = (1+O(n^{-2\bb}))\prod_{ k=l}^{ n-1}{m_{ k+1}+m_{ k}-1\choose m_{
k}-1}\( { 1\over 2} \)^{m_{ k+1}+m_{ k} }.\nn
\end{eqnarray} Our analysis also shows that for some $\de_{3, l}\rar 0$
\begin{equation}
         \inf_{\stackrel{ m_{ l}} {m_{l}\stackrel{l}{\sim}\NNN_{ l}}
}\left(\sum_{\stackrel{m_{ 2}, \ldots, m_{ l-1}} {m_{j}\stackrel{j}{\sim}
\NNN_{ j}} }\prod_{ k=2}^{ l-1}{m_{ k+1}+m_{ k}-1\choose m_{ k}-1}\( {
1\over 2}
\)^{m_{ k+1}+m_{ k} }\right)\geq ( ( l-1)!)^{ -3a-\de_{3, l}},\label{44.8y}
\end{equation} and since
\begin{eqnarray} &&
\sum_{\stackrel{m_{ 2}, \ldots, m_{ n}} {m_{j}\stackrel{j}{\sim}\NNN_{ j}} }
\prod_{\ell= 2}^{ n-1}{m_{ \ell+1}+m_{ \ell}-1\choose m_{
\ell}-1}\( { 1\over 2} \)^{m_{ \ell+1}+m_{ \ell} }\label{44.8z}\\
&&\geq\sum_{\stackrel{m_{ l}, \ldots, m_{ n}} {m_{j}\stackrel{j}{\sim}\NNN_{
j}} }\prod_{ k=l}^{ n-1}{m_{ k+1}+m_{ k}-1\choose m_{ k}-1}\( { 1\over 2}
\)^{m_{ k+1}+m_{ k} }\nn\\ &&
         \inf_{\stackrel{ m_{ l}} {m_{\ell}\stackrel{\ell}{\sim}\NNN_{\ell}}
}\left(\sum_{\stackrel{m_{ 2}, \ldots, m_{ l-1}} {m_{j}\stackrel{j}{\sim}\NNN
_{ j}} }
\prod_{ k=2}^{ l-1}{m_{ k+1}+m_{ k}-1\choose m_{ k}-1}\( { 1\over 2}
\)^{m_{ k+1}+m_{ k} }\right),\nonumber
\end{eqnarray} where we used the fact that for $C( i,j),D( i,j)$ 
non-negative, we
have
$$\sum_{ i,j,k}C( i,j)D( j,k) =\sum_{j}\sum_{ i}C( i,j)\sum_{k}D(  j,k)
\geq \(\sum_{ i,j}C( i,j)\) \inf_{ j}\sum_{k}D(  j,k),$$ we see from 
(\ref{44.8})
and (\ref{44.6f}) that uniformly in
$z\in \partial D(x,r_{n,l})_{ n^{ 4}}$
\begin{equation}
         \sum_{\stackrel{m_{ l}, \ldots, m_{ n}} {m_{k}\stackrel{k}{\sim}\NNN_{
k}} }\PPP^{ z}
\( N_{n,l,m_{ l},k}^x =m_{k},\; k=l+1,\ldots,n\,;\,\Om_{ x,n,l,m_{
l}}^{l-1,\ldots,n}\)
\leq C\log n\, Q_n  \,\,( l!)^{3a+\de_{3,l }}.\label{44.9s}
\end{equation} As in (\ref{pr.4}) we have that uniformly in $n$ and
$x\in U_{ n}$
\begin{equation}
\PPP \(T_{D(x,r'_{n,l})}<T_{ D(0,K_n)^{ c}}\,;
\,T_{D(x,r_{n,l}')}=T_{\partial D(x,r_{n,l})_{ n^{ 4}} }\) \leq c'/l\log
n\label{pr.4rep}
\end{equation} so that by readjusting $\de_{3,l }$
\begin{equation}
         \sum_{\stackrel{m_{ l}, \ldots, m_{ n}} {m_{k}\stackrel{k}{\sim}\NNN_{
k}} }\PPP
\( N_{n,l,m_{ l},k}^x =m_{k},\; k=l+1,\ldots,n\,;\,\Om_{
x,n,l,m_{l}}^{l-1,\ldots,n}\)
\leq C\, Q_n  \,\,( l!)^{3a+\de_{3,l }},\label{44.9q}
\end{equation} and (\ref{44.9qd}) follows.\qed

\begin{lemma}\label{cond-corm} For some $C<\ff$ and $\de_{3, l}\rar 0$
\bea &&
         \sum_{\stackrel{m_{ 2}, \ldots, m_{ l}} {m_{k}\stackrel{k}{\sim}\NNN_{
k}} }\PPP
\( N_{n,k}^x =m_{k},\; k=2,\ldots,l\,;\,\Om_{ x,n,1,1}^{1,\ldots,l }\)
\leq C\,( l!)^{-3a+\de_{3,l }}.\label{44.9qe}
\eea
\end{lemma}

\noindent{\bf Proof of Lemma \ref{cond-corm}:} As before, uniformly in
$m_{k}\stackrel{k}{\sim}\NNN_{k},\; k=2,3,\ldots,l$ and $z\in
\partial D(x,r_{n,1})_{ n^{ 4}}$
\begin{eqnarray}
         &&
         \PPP^{ z} \( N_{n,1,1,k}^x =m_{k},\; k=2,\ldots,l\,;\,\Om_{
x,n,1,1}^{1,\ldots,l }\)
\label{44.8m}\\
         &&
         = (1+O(n^{-2\bb}))\prod_{ k=2}^{ l-1}{m_{ k+1}+m_{ k}-1\choose m_{
k}-1}\( { 1\over 2} \)^{m_{ k+1}+m_{ k} }.\nn
\end{eqnarray} Using (\ref{44.5}) as before, we obtain (\ref{44.9qe}).
\qed

\section{Second moment estimates}\label{52momest}

We begin by defining the
$\si$-algebra $\GG_{n, l}^x$ of excursions from
$D(x,r_{n, l-1})^{ c}$ to $D(x,r'_{n, l})$. To this end, fix
$x \in \Z^2$, let $\taup_0=0$ and for $i=1,2,\ldots$ define
\begin{eqnarray*}
\tau_{i} & = & \inf \{ k \geq \taup_{i-1}  :\; X_k \in D(x,r'_{n, l})\}
\,, \\
\taup_{i} & = & \inf \{ k \ge \tau_{i} :\; X_k \in  D(x,r_{n, l-1})^{ c}\}.
\end{eqnarray*} Then $\GG_{n, l}^x$ is the $\sigma$-algebra generated by
the excursions
$\{ e^{(j)}, j =1,\ldots \}$, where
$e^{(j)}=\{ X_{k} : \taup_{j-1}\leq k \leq \tau_{j} \}$ is the $j$-th excursion
from
$D(x,r_{ n,l-1})^{ c}$ to $D(x,r'_{n, l})$ (so for $j=1$ we do begin at
$t=0$).

The following Lemma is proved in the next section. Recall that for any
$\si$-algebra $\GG$ and event $B\in\GG$, we have
         $\P( A \cap B\,|\, \GG)=\P( A\,|\, \GG)1_{ \{ B\}}$.

	\begin{lemma}[Decoupling Lemma]\label{cond-cor} Let
	\[\Gamma_{n,l}^y =\{ N^y_{n,i} = m_i ; i=l+1,\ldots,n \}\cap
\Om_{ x,n,l,m_{
	l}}^{l-1,\cdots,n }.\]  Then, uniformly over all $l\leq n ,\,
	m_l \stackrel{l}{\sim} \NNN_l$,  $\{m_i: i=l,\ldots,n\}$,
	$y \in U_{ n}$,
	\begin{eqnarray} &&
	\PPP ( \Gamma_{ n,l}^y\, ,\,N_{n,l}^y=m_l\,
	|\, {\mathcal G}_{ n,l}^y)\label{new1.3ez}\\ &&
	  = (1+O(n^{-1/2}))  \PPP (\Gamma_{ n,l}^y
	\, |\,N_{n,l}^y=m_l)1_{ \{ N_{n,l}^y=m_l\}} \nonumber
	\end{eqnarray}
	\end{lemma}

	\begin{remark}{\rm The intuition behind the Decoupling Lemma 
is that what
	happens `deep inside' $D(y,r'_{n, l})$, e.g., $\Gamma_{ n,l}^y $, is
	`almost' independent of what happens outside $D(y,r'_{n, l})$, i.e.,
	${\mathcal G}_{
	n,l}^y$.}
	\end{remark}

	\noindent{\bf Proof of (\ref{pm3.5}):} Recall that $\NNN_k=3a k^2
	\log k$ and
	that we write
	$m\stackrel{k}{\sim} \NNN_{k}$ if $m=1$ for $k<k_{ 0}$ and
	$|m-\NNN_k| \leq k$ for $k\geq k_{ 0}$. Relying upon the first
	moment estimates
	and Lemma \ref{cond-cor}, we next prove the second moment estimates
\req{pm3.5}. Take $x,y \in U_{n}$ with $l( x,y)=l-1$. Thus $2 r_{n,l-1}+2
\leq |x-y|< 2 r_{n,l-2}+2$ for some $2 \leq l \leq n$. Since
$r_{n,l-3} -  r_{n,l-2} \gg  2 r_{n,l-1}$, it is easy to see that
$\partial D(y,r_{n,l-1})_{ n^{ 4}} \cap \partial D(x,r_{n,k})_{ n^{ 4}}
=\emptyset$ for all $k
\neq l-2$.
         Replacing hereafter $l$ by $l
\wedge (n-3)$, it follows that for
$k\neq l-1,l-2$,
         the events $\{ N_{n,k}^x \stackrel{k}{\sim} \NNN_{k} \}$ are measurable
with respect to  the
$\sigma$-algebra $\GG^y_{n,l}$.

We write
\[\Gamma_{n,l}^y( m_{ l},\ldots, m_{ n}) =\{ N^y_{n,i} = m_i ; i=l+1,\ldots,n
\}\cap
\Om_{ x,n,l,m_{ l}}^{l-1,\cdots,n }\] to emphasize the dependence on
$ m_{ l},\ldots, m_{ n}$. With $J_l:=\{ l+1,\ldots,n \}$ set
\[\wt{\Ga}_{ n}^{ y}(J_l,m_{ l} )=\bigcup_{m_{k }
\stackrel{k}{\sim}
\NNN_{k}\,;\,k
\in J_{l}}
\Gamma_{n,l}^y( m_{ l},\ldots, m_{ n})\] Similarly, with
$M_{l-3 } :=\{ 2,\ldots,l-3 \}$ set
         \[\wt{\Ga}^{ x}_{ n}(M_{l-3 })=\{ N_{n,k}^x \stackrel{k}{\sim}
\NNN_{k},\; k \in M_{l-3 }
\}\cap
\Om^{1,\ldots , l-3}_{ x,n,1,1},\]
         and with
$I_l:=\{ 2,\ldots,l-3,l,\ldots,n \}$ set
\[\bar{\Ga}_{ n}^{ x}(I_l )=\bigcup_{m'_{l} \stackrel{l}{\sim}
\NNN_{l}} \wt{\Ga}_{ n}^{ x}(J_l,m'_{ l} )\cap \{ N^x_{n,l} =m'_{l}\}\cap
\wt{\Ga}_{ n}^{  x}(M_{l-3 }).\] Using the previous paragraph we can check that
$\bar{\Ga}_{ n}^{ x}(I_l )\in {\mathcal G}^y_{ n,l}$. Note that
\begin{eqnarray}  &&
\{ x, y \, \mbox{ are $n$-successful} \}\label{jstar}\\ &&\hspace{ .5in}
\subseteq \bigcup_{m_l\stackrel{l}{\sim} \NNN_{l} }
\lc\bar{\Ga}_{ n}^{ x}(I_l )\bigcap \wt{\Ga}_{ n}^{ y}(J_l,m_{ l} )\bigcap
\{ N^y_{n,l} = m_{l}\}\rc.\nn
\end{eqnarray} Applying \req{new1.3ez}, we have that for some universal
constant
$C_3<\infty$,
\begin{eqnarray} &&
\PPP\(x\, \mbox{ and } \, y \, \mbox { are $n$-successful}
\)\label{52mom}\\ &&
\leq    \sum_{m_l\stackrel{l}{\sim} n_{l}}
\E \left[ \PPP ( \wt{\Ga}_{ n}^{ y}(J_l,m_{ l}  )\,,N^y_{n,l} = m_{l}
\,
\big|
\,
         {\mathcal G}^y_{ n,l} ) \,;\bar{\Ga}_{ n}^{ x}(I_l )
\right]\nonumber\\ &&\leq C_3
\PPP (\bar{\Ga}_{ n}^{ x}(I_l ), N^y_{n,l} = m_{l} )
\sum_{m_l\stackrel{l}{\sim} n_{l}}
          \PPP ( \wt{\Ga}_{ n}^{ y}(J_l,m_{ l}  )\,
\big| \, N^y_{n,l} = m_{l} )\nn\\ &&\leq C_3
\PPP (\bar{\Ga}_{ n}^{ x}(I_l ))
\sum_{m_l\stackrel{l}{\sim} n_{l}}
          \PPP ( \wt{\Ga}_{ n}^{ y}(J_l,m_{ l}  )\,
\big| \, N^y_{n,l} = m_{l} ).\nn
\end{eqnarray} Using (\ref{44.9qd}), for some universal constant
$C_5<\infty$,
\begin{eqnarray}\qquad
\sum_{ m_{l}\stackrel{l}{\sim} n_{l}}\PPP (
         \wt{\Ga}_{ n}^{ y}(J_l,m_{ l}  )\, \big|
           \, N^y_{n,l} = m_{l} )
         \leq C_5\, Q_n  \,\,( l!)^{3a+\de_{3,l }}. \label{244.9}
\end{eqnarray}

Noting that  $\wt{\Ga}_{ n}^{ x}(M_{l-3 })\in  {\mathcal G}^x_{  n,l}$,
(\ref{new1.3ez}) then shows that
\begin{eqnarray} &&
\PPP \( \bar{\Ga}_{ n}^{ x}(I_l )\)\label{244.7}\\
           &&\leq
\sum_{m_l\stackrel{l}{\sim} \NNN_{l}}
\E \left[ \PPP ( \wt{\Ga}_{ n}^{ x}(J_l,m_{l} ),\,N^x_{n,l} = m_{l} \,
\big|
\,
         {\mathcal G}^x_{ n,l} ) \,;
\wt{\Ga}_{ n}^{ x}(M_{l-3 }) \right]\nn\\
           &&\leq C_6 \PPP \( \wt{\Ga}_{ n}^{ x}(M_{l-3 } ), N^x_{n,l} 
= m_{l} \)
\sum_{ m_{l}\stackrel{l}{\sim} \NNN_{l}}
\PPP ( \wt{\Ga}_{ n}^{ x}(J_l,m_{l}  )\, \big|
           \, N^x_{n,l} = m_{l} )\nn\\
           &&\leq C_6 \PPP \( \wt{\Ga}_{ n}^{ x}(M_{l-3 } )\)
\sum_{ m_{l}\stackrel{l}{\sim} \NNN_{l}}
\PPP ( \wt{\Ga}_{ n}^{ x}(J_l,m_{l}  )\, \big|
           \, N^x_{n,l} = m_{l} ).\nn
\end{eqnarray} Using  (\ref{44.9qe}) and  (\ref{244.9}) we get that, for some
$\de_{4,l}\rar 0$
\begin{equation}
\PPP \(\bar{\Ga}_{ n}^{ x}(I_l ) \)
\leq C_7 l^{ 15}\,  ( l!)^{\de_{4,l }}\,Q_n.\label{244.10}
\end{equation} Putting (\ref{52mom}), (\ref{244.9}) and (\ref{244.10})
together and adjusting $C$ and $\de'_{ l-1}$ proves (\ref{pm3.5}) for
$l( x,y)=l-1$.
\qed

\section{Approximate decoupling}\label{sec-appdec}

The goal of this section is to prove the Decoupling Lemma, Lemma
\ref{cond-cor}. Since what happens `deep inside' $D(y,r'_{n, l})$, e.g.,
$\Gamma_{ n,l}^y $, depends on what happens outside $D(y,r'_{n, l})$, i.e., on
${\mathcal G}_{ n,l}^y$, only through the initial and end points of the
excursions from
$D(x,r'_{ n,l})$ to $D(x,r_{n, l-1})^{ c}$, we begin by studying the dependence
on these initial and end points.

Consider a random path beginning at $z
\in  \partial D(0,r_{n, l})_{ n^{ 4}}$. We will show that for $l$ 
large, a certain
$\si$-algebra of excursions of the path from $ D(0,r'_{n, l+1})$ to
$D(0,r_{ n, l})^{ c}$ prior to $T_{D(0,r_{ n,l-1})^{ c} }$, is almost 
independent
of the  choice of initial point
$z
\in  \partial D(0,r_{ n,l})_{ n^{ 4}}$ and  final point $w \in
\partial D(0,r_{n, l-1})_{n^{ 4}}$.  Let
$\tau_0=0$ and for $i=0,1,\ldots$ define
\begin{eqnarray*}
\tau_{2i+1} & = & \inf \{ k \geq \tau_{2i}  :\; X_k \in  D(0,r'_{n,  l+1})
\cup  D(0,r_{ n,l-1})^{ c} \} \\
\tau_{2i+2} & = &
           \inf \{ k \ge \tau_{2i+1} :\; X_k \in  D(0,r_{ n,l})^{ c}\}\,.
\end{eqnarray*} Abbreviating $\bar{\tau}=T_{ D(0,r_{n, l-1})^{ c}}$ note that
$\bar{\tau}=\tau_{2I+1}$ for some (unique) non-negative integer
$I$. As usual, $\FF_{j}$ will denote the $\si-$algebra generated by
$\{X_{ l},\,l=0,1,\ldots, j \}$, and for any stopping time $\tau$,
$\FF_{\tau}$ will denote the collection of events $A$ such that
$A\cap \{ \tau=j\}\in \FF_{j}$ for all $j$.

Let $\HH_{n, l}$ denote the $\si$-algebra generated by the excursions of the
path from
$ D(0,r'_{n,  l+1})$ to  $D(0,r_{n, l})^{ c}$, prior to $T_{ D(0,r_{ 
n,l-1})^{ c} }$.
Then
$\HH_{ n,l}$ is the $\sigma$-algebra generated by the excursions
$\{ v^{(j)}, j =1,\ldots, I \}$, where
$v^{(j)}=\{ X_{k} : \tau_{2j-1}\leq k \leq \tau_{2j} \}$ is the
$j$-th excursion  from $ D(0,r'_{n, l+1})$ to  $D(0,r_{n, l})^{ c}$.

\begin{lemma}\label{lemprobends}  Uniformly in   $l,n$, $z,z'\in
\partial D(0,r_{n, l})_{ n^{ 4}}$,
$w\in  \partial D(0,r_{ n,l-1})_{ n^{4}}$, and $B_{ n} \in \HH_{n, l}$,
\bea &&
\PPP^z (B_{ n}\cap \Om^{ l-1,l,l+1}_{ 0,n,l,1} \,\big |\,X_{T_{ 
D(0,r_{n, l-1})^{
c} }}=w)\label{m1.5euj}\\ && = (1 +O(n^{ -3}))
\PPP^{ z} (B_{ n}\cap
\Om^{ l-1,l,l+1}_{ 0,n,l,1}),
\nn
\eea and
\begin{equation}
\PPP^{ z} (B_{ n}\cap \Om^{ l-1,l,l+1}_{ 0,n,l,1}) = (1 +O(n^{ -3}))
\PPP^{ z'} (B_{ n}\cap
\Om^{ l-1,l,l+1}_{ 0,n,l,1})
\,.
\label{m1.5evj}
\end{equation}
\end{lemma}

\noindent{\bf Proof of Lemma \ref{lemprobends}:} Fixing $z \in
\partial D(0,r_{n, l})_{ n^{ 4}}$ it suffices to consider
$B_{ n} \in \HH_{n, l}$ for which $\PPP^z(B_{ n})>0$. Fix such a set
$B_{ n}$ and a point
$w\in
\partial D(0,r_{n, l-1})_{ n^{ 4}}$. Using the notation introduced 
right before
the statement of our Lemma, for any $i
\geq 1$, we can write
\bea &&\{ B_{ n}\cap \Om^{ l-1,l,l+1}_{ 0,n,l,1}, I=i \}\nn\\ &&=\{ B_{n,
i}\cap A_{ i},
\tau_{2i}<\bar{\tau}\}\cap(\{I=0\,, X_{\bar{\tau} }\in  \partial D(0,r_{n,
l-1})_{ n^{ 4}}\}\circ\th_{\tau_{2i}})\nn\eea
         for some $B_{n, i}\in\FF_{ \tau_{2i}}$, where
\[A_{ i}=\lc X_{\tau_{2j-1} }\in
\partial D(0,r_{n, l+1})_{ n^{ 4}}\,, X_{\tau_{2j} }\in  \partial 
D(0,r_{ n,l})_{
n^{ 4}}\,,
\,\forall\,j\leq i\rc\in\FF_{
\tau_{2i}}\]
         so by the strong Markov property at $\tau_{2i}$,
\bea &&
\E^z [ X_{\bar\tau}=w ; B_{ n}\cap \Om^{ l-1,l,l+1}_{ 0,n,l-1,1}, I=i 
]\nn\\ && =
\E^z \left[ \E^{X_{\tau_{2i}}} (X_{\bar\tau}=w,\, I=0) ; B_{n, i}\cap A_{i},
\tau_{2i} < {\bar\tau} \right] \,,
\nn\eea and
\bea &&
\PPP^z \(B_{ n}\cap \Om^{ l-1,l,l+1}_{ 0,n,l,1}, I=i \)\nn\\ && =
\E^z
\left[
\E^{X_{\tau_{2i}}} (I=0\,, X_{\bar{\tau} }\in  \partial D(0,r_{ 
n,l-1})_{ n^{ 4}})
; B_{n, i}\cap A_{ i}, \tau_{2i} < {\bar\tau} \right]
           \;.
\nn\eea Consequently, for all $i \geq 1$,
\begin{eqnarray} &&
\E^z [ X_{\bar\tau}=w;  B_{ n}\cap \Om^{ l-1,l,l+1}_{ 0,n,l,1}, I=i
]\label{ofer-star2}\\ &&
         \geq \PPP^z \(B_{ n}\cap \Om^{ l-1,l,l+1}_{ 0,l,1}, I=i 
\)\nn\\ &&\hspace{
1in} \times
\inf_{x\in \partial D(0,r_{n, l})_{ n^{ 4}}}
\frac{\E^{x} \( X_{\bar\tau}=w;\, I=0 \)} {\E^x \(I=0\,,  X_{\bar{\tau} }\in
\partial D(0,r_{n, l-1})_{ n^{ 4}}\)}
\,. \nonumber
\end{eqnarray} Necessarily $\PPP^z(B_{ n}|I=0)\in
\{0,1\}$ and is independent of $z$ for any $B_{ n} \in \HH_{ n,l}$, implying
that
\req{ofer-star2} applies for
$i=0$ as well. By
\req{p1.5h}, \req{p1.5}, \req{p1.2bfa}  and \req{p1.2bd} there exists
$c<\infty$ such that for any $z,x\in \partial D(0,r_{n, l})_{ n^{ 4}}$ and
$w \in
\partial D(0,r_{n, l-1})_{ n^{ 4}}$,
\[
\frac{\E^{x} \( X_{\bar\tau}=w;\, I=0 \)} {\E^x \(I=0\,,  X_{\bar{\tau} }\in
\partial D(0,r_{ n,l-1})_{ n^{ 4}}\)} \geq (1-c n^{ -3})  H_{D(0,r_{ n,l-1})^{
c}}(z,w)\,.
\] Hence, summing \req{ofer-star2} over $I=0,1,\ldots$, we get that
\bea &&
           \E^z \left[ X_{\bar\tau}=w, B_{ n}\cap \Om^{ l-1,l,l+1}_{ 0,n,l,1}
\right]\label{declb}\\ &&\hspace{.8in}
\geq (1 -c n^{ -3})\PPP^z (B_{ n}\cap \Om^{ l,l+1}_{ 0,n,l,1}) H_{ D(0,r_{n,
l-1})^{ c}}(z,w)\,.\nn
\eea A similar argument shows that
\bea &&
\E^z \left[  X_{\bar\tau}=w, B_{ n}\cap \Om^{ l-1,l,l+1}_{ 0,n,l,1}
\right]
\label{decub}\\ &&\hspace{.8in}\leq (1+c n^{ -3}) \PPP^z \(B_{ n}
\cap
\Om^{ l-1,l,l+1}_{ 0,n,l,1}\) H_{ D(0,r_{n, l-1})^{ c}}(z,w)\,,
\nn
\eea and we thus obtain (\ref{m1.5euj}).

By the strong Markov property at $\tau_1$, for any $z\in   \partial D(0,r_{n,
l})_{ n^{ 4}}$,
\bas &&
\PPP^{ z} (B_{ n}\cap \Om^{ l,l+1}_{ 0,n,l,1})= \PPP^{z}(B_{ n}\cap
\Om^{ l-1,l,l+1}_{ 0,l,1}, I=0)
\\ &&\hspace{.2in}+\sum_{ x\in  \partial D(0,r_{n, l+1})_{ n^{ 4}}}H_{
D(0,r'_{ l+1})\cup   D(0,r_{n, l-1})^{ c}}(z,x)\PPP^{x}\(B_{ n}\cap
\Om^{ l-1,l,l+1}_{ 0,n,l,1}\)
\eas The term involving $\{B_{ n}\cap \Om^{ l-1,l,l+1}_{ 0,n,l,1} , I=0\}$ is
dealt with by (\ref{p1.2b}) and (\ref{m1.5evj}) follows by
\req{p1.5b}.
\qed

Building upon Lemma \ref{lemprobends} we quantify the independence
between the
$\si$-algebra $\GG_{ l}^x$ of excursions from
$D(x,r_{ n,l-1})^{ c}$ to $D(x,r'_{n, l})$ and the
$\sigma$-algebra
$\HH_{n, l}^x (m)$ of excursions from $D(x,r'_{ n, l+1})$ to
$D(x,r_{n, l})^{ c}$ during  the first $m$ excursions from
         $D(x,r'_{n, l})$ to $D(x,r_{ n,l-1})^{ c}$. To this end, fix
$x \in \Z^2$, let $\taup_0=0$ and for $i=1,2,\ldots$ define
\begin{eqnarray*}
\tau_{i} & = & \inf \{ k \geq \taup_{i-1}  :\; X_k \in D(x,r'_{ n,l}) \}
\,, \\
\taup_{i} & = & \inf \{ k \ge \tau_{i} :\; X_k \in  D(x,r_{ n,l-1})^{ c}
\}.
\end{eqnarray*} Then $\GG_{ l}^x$ is the $\sigma$-algebra generated by the
excursions
$\{ e^{(j)}, j =1,\ldots \}$, where
$e^{(j)}=\{ X_{k} : \taup_{j-1}\leq k \leq \tau_{j} \}$ is the $j$-th excursion
from
$D(x,r_{ l-1})^{ c}$ to $D(x,r'_{n, l})$ (so for $j=1$ we do begin at
$t=0$).

We denote by
$\HH_{ n,l}^x (m)$ the $\sigma$-algebra generated by all excursions from
$D(x,r'_{ n,l+1})$ to
$D(x,r_{n, l})^{ c}$ from time $\tau_1$ until time $\taup_{m}$. In more detail,
for each
$j=1,2,\ldots,m$ let $\sip_{j,0}=\tau_{j}$ and for $i=1,\ldots$ define
\begin{eqnarray*}
\sii_{j,i} & = & \inf \{ k \geq \sip_{j,i-1}  :\; X_k \in
         D(x,r'_{ n, l+1})  \}\,, \\
\sip_{j,i} & = &
           \inf \{ k \ge \sii_{j,i} :\; X_k \in  D(x,r_{n, l})^{ c}\}\,.
\end{eqnarray*} Let $v_{j,i}=\{ X_{k} : \sii_{j,i}\leq k \leq
\sip_{j,i} \}$ and
$Z^{j}=\sup \{ i \geq 0 \,:\,  \sip_{j,i}<\taup_{j}\}$. Then, $\HH_{ n,l}^x(m)$
is the
         $\sigma$-algebra generated by the intersection of the
$\si$-algebras
$\HH_{n, l,j}^x =
\sigma(v_{j,i}, i=1,\ldots,Z^{j})$ of the excursions between times
$\tau_j$ and
$\taup_j$, for $j=1,\ldots,m$.

\begin{lemma}\label{cond-ind} There exists $C<\infty$ such that uniformly
over all
$m \leq (n\log n)^{2}$, $l$,
$x \in \Z^2$ and $y_0,y_1 \in \Z^2 \setminus D(y,r'_{n,l})$, and
$H
\in \HH_{n, l}^x (m)$,
\begin{eqnarray} && (1-C m n^{ -3} )  \PPP^{y_1} (H\cap
\Om^{l-1,l,l+1}_{ x,n,l,m})
\leq
\PPP^{y_0} (H\cap \Om^{l-1,l,l+1}_{ x,n,l,m} \, |\, \GG_{ l}^x
)\label{new1.2e}\\ &&
\hspace{ 2in}\leq (1+C m n^{ -3}  )  \PPP^{y_1} (H\cap
\Om^{l-1,l,l+1}_{ x,n,l,m})
          \,. \nonumber
\end{eqnarray}
\end{lemma}

\noindent{\bf Proof of Lemma \ref{cond-ind}:} Applying the Monotone Class
Theorem to the algebra of their finite disjoint unions,  it suffices to prove
\req{new1.2e} for the generators of the
$\si$-algebra $\HH_{ n,l}^x (m)$ of the form $H=H_1 \cap H_2
\cap
\cdots
\cap H_m$, with $H_j \in \HH^x_{n,l,j}$ for $j=1,\ldots,m$. Conditioned upon
$\GG_{ l}^x$  the events $H_j$ are independent. Further, each
$H_j$ then has the conditional law of an event $B_j$ in the
$\si$-algebra
$\HH_{ n,l}$ of Lemma
\ref{lemprobends}, for some random $z_j =X_{\tau_j}-x \in
\partial D(0,r_{ n,l})_{ n^{ 4}}$ and
$w_j =X_{\taup_j}-x \in \partial D(0,r_{n, l-1})_{ n^{ 4}}$, both measurable on
$\GG_{ l}^x$. By our conditions, the uniform estimates (\ref{m1.5euj}) and
(\ref{m1.5evj}) yield that for any fixed $z'\in
\partial D(0,r_{ n,l})_{ n^{ 4}}$,
\beqn{ofer-star1} &&
\PPP^{y_0} (H\cap \Om^{l-1,l,l+1}_{ x,n,l,m} \, |\, \GG_{ l}^x )\\
&&=\PPP^{y_0} (\cap _{ j=1}^{ m}(H_j\cap \Om^{l-1,l,l+1}_{ x,n,l,1} )
\, |\,
\GG_{ l}^x )\nn\\ && = \prod_{j=1}^m \PPP^{z_j}(B_j \cap
\Om^{l-1,l,l+1}_{ x,n,l,1} \, | \, X_{T_{D(0,r_{ l})^{ c}}} =w_j ) \nn
\\ && =
\prod_{j=1}^{m} \,(1 + O( n^{ -3}))\PPP^{z_j}(B_j\cap
\Om^{l-1,l,l+1}_{ 0,n,l,1} )\nn \\ && = (1 +O(n^{ -3} ) )^{m}
\prod_{j=1}^{m} \,\PPP^{z'}(B_j\cap \Om^{l-1,l,l+1}_{ 0,n,l,1} )\,.\nn
\eeqn Since $m \leq (n\log n)^{2}$ and the right-hand side of
\req{ofer-star1} neither depends on
$y_0 \in \Z^2$ nor on the extra information in $\GG_{ l}^x$, we get
\req{new1.2e}.
\qed

\begin{corollary}\label{cond-cora} Let $\Gamma_{n,l}^y =\{ N^y_{n,i} = m_i ;
i=l+1,\ldots,n \}\cap
\Om^{l-1,\cdots,n}_{ y,n,l,m_{l}}$. Then, uniformly over all $n
\geq l ,\,$
$m_l \stackrel{l}{\sim} \NNN_l$, $\{m_i: i=l,l+2,\ldots,n\}$,
$y \in U_{ n}$ and $x_0,x_1 \in \Z^2 \setminus D(y,r'_{n,l})$,
\begin{eqnarray} &&
\PPP^{x_0} ( \Gamma_{n,l}^y\,,\, N_{n,l}^y=m_l\,|\, {\mathcal G}_l^y)
\label{new1.3e}\\ &&= (1+O(n^{-1}  \log n ))  \PPP^{x_1} (\Gamma_{n,l}^y
\, |\,N_{n,l}^y=m_l){\bf 1}_{\{N^y_{n,l}=m_l\}}\nn
\end{eqnarray}
\end{corollary}

\noindent{\bf Proof of Corollary \ref{cond-cora}:} For
$j=1,2,\ldots$ and
$i=l+1,\ldots,n$, let $Z_i^j$ denote the number of excursions from
$ D(y,r_{n,i-1})^{ c}$ to $D(y,r'_{n,i})$ by the random walk during the time
interval
$[\tau_j,\taup_j]$.  Clearly, the event
$$ H =\{ \sum_{j=1}^{m_l} Z_i^j = m_i : i=l+1,\ldots,n \}\cap
\Om^{l+1,\ldots,n}_{ y,n,l+1,m_{l+1}}
$$ belongs to the $\si$-algebra $\HH_{ n,l}^y (m_l)$ of Lemma
\ref{cond-ind}. It is easy to verify that starting at any $x_{ 0} \notin
D(y,r'_{n,l}) $,  when the event
$\{N^y_{n,l}=m_l\} \in \GG^y_l$ occurs, it implies that
$N^y_{n,i}=\sum_{j=1}^{m_l} Z_i^j$ for $i=l+1,\ldots,n$. Thus,
\begin{eqnarray} &&
\PPP^{x_0} ( \Gamma_{n,l}^y \, |{\mathcal G}_l^y) {\bf
1}_{\{N^y_{n,l}=m_l\}} =
\PPP^{x_0} ( H\cap \Om^{l-1,l,l+1}_{ y,n,l,m_{ l}}\, |{\mathcal G}_l^y) {\bf
1}_{\{N^y_{n,l}=m_l\}}\,. \label{hole0}
\end{eqnarray} With $m_l/(n^2 \log n)$ bounded above, by
\req{new1.2e} we have, uniformly in
$y \in \Z^2$ and $x_0,x_1 \in \Z^2 \setminus D(y,r'_{n,l})$,
\begin{equation}
\label{hole1}
\PPP^{x_0} ( H\cap \Om^{l-1,l,l+1}_{ y,n,l,m_{ l}}\, |{\mathcal G}_l^y)  =
         (1+O( n^{-1}  \log n ))
\PPP^{x_1} ( H\cap \Om^{l-1,l,l+1}_{ y,n,l,m_{ l}})\,.
\end{equation} Hence,
\begin{eqnarray} &&
\PPP^{x_0} ( \Gamma_{n,l}^y\, |{\mathcal G}_l^y)  {\bf
1}_{\{N^y_{n,l}=m_l\}}\label{hole2a}\\ && =(1+O( n^{-1}
\log n ))
\PPP^{x_1} ( H\cap \Om^{l-1,l,l+1}_{ y,n,l,m_{ l}}) {\bf
1}_{\{N^y_{n,l}=m_l\}}\,.
\nonumber
\end{eqnarray}
       Setting $x_0=x_1$ and taking expectations with respect to
$\PPP^{x_0}$, one has
\begin{eqnarray} &&
\PPP^{x_1} ( \Gamma_{n,l}^y \, | N^y_{n,l}=m_l)= (1+O(n^{-1}
\log n ))
\PPP^{x_1} ( H\cap \Om^{l-1,l,l+1}_{ y,n,l,m_{ l}}).
\label{hole2}
\end{eqnarray} Hence,
\begin{eqnarray}&&
\PPP^{x_1} ( \Gamma_{ n,l}^y\, | N^y_{n,l}=m_l){\bf
1}_{\{N^y_{n,l}=m_l\}}\label{hole2f}\\&&=(1+O(n^{-1}
\log n))
\PPP^{x_1} ( H \cap \Om^{l-1,l,l+1}_{ y,n,l,m_{ l}}){\bf 1}_{\{N^y_{n,l}=m_l\}}
\nn \\&&=(1+O(n^{-1}
\log n)) \PPP^{x_0} ( \Gamma_{ n,l}^y\, |{\mathcal G}_{n,l}^y)  {\bf
1}_{\{N^y_{n,l}=m_l\}} \nonumber
\end{eqnarray} where we used (\ref{hole2a}) for the last equality.
         Using that
${\{N^y_{n,l}=m_l\}}\in {\mathcal G}_{n,l}^y\,$, this is
\req{new1.3e}.
\qed

\section{Appendix}

   Let $q_0(x)={\bf 1}_{\{0\}}(x)$ and for $n\geq 1$ let
$$q_n(x)=\frac{1}{2\pi n} e^{-|x|^2/2n}.$$

\begin{proposition}\label{AP1} Suppose $X_n$ is a strongly aperiodic
symmetric random walk in $\Z^2$ with the covariance matrix of $X_1$ equal
to the identity and with
$3+2\beta$ moments. Then there exists $c_1$ such that
$$\sup_{x\in \Z^2} |p_n(x)-q_n(x)|\leq c n^{-\frac32-\beta}, \qquad n\geq
1.$$
\end{proposition}

\proof Let $\phi$ be the characteristic function for $X_1$. Since $X_1$ is
symmetric, the third moments are zero, that is, if
$X_1=(X_1^{(1)},X^{(2)}_1)$, $i_1, i_2\geq 0$, $i_1+i_2=3$, then
$\E[(X_1^{(1)})^{i_1}(X_1^{(2)})^{i_2}]=0$. So by a Taylor expansion,
$$\phi(\al/\sqrt n)=1-\frac{|\al|^2}{n} +E_1(\al,n),$$ where
$$|E_1(\al,n)|\leq c_2 (|\al|/\sqrt n)^{3+2\beta},$$  provided $\al\in
[-\pi,\pi]^2$. Similarly
$$e^{-|\al|^2/n}=1-\frac{|\al|^2}{n} +E_2(\al,n),$$ where the error term
$E_2(\al,n)$ has the same bound.  We now follow Proposition 3.1 of
\cite{locur}, using the above estimate for $E_i(\al,n)$, $i=1,2$, in 
place of the
one in that paper.
\qed

\begin{proposition}\label{AP2} Let $X_n$ be as above and
$$a(x)=\sum_{n=0}^\infty [p_n(0)-p_n(x)].$$ Then for $x\ne 0$
\begin{equation}\label{AE1} a(x)=\frac{2}{\pi} \log |x| + k + o(1/|x|)
\end{equation} where $k$ is a constant depending on $p_1$ but not $x$.
\end{proposition}

\proof We write
\begin{align} a(x)&=\sum_n[p_n(0)-q_n(0)] +\sum_n
[q_n(0)-q_n(x)]\label{AE2}\\ &~~~+\sum_n [q_n(x)-p_n(x)]\nn\\ &=
I_1+I_2+I_3.\nn
\end{align} Since $n^{-\frac32-\beta}$ is summable, $I_1$ is a constant not
depending on $x$. $I_2$ has the form given in the right hand side of
(\ref{AE1}); see the proof of Theorem 1.6.2 in \cite{L1}. So it remains to show
$I_3=o(1/|x|)$.

We write
\begin{align} I_3&=\sum_{n\leq N} p_n(x) +\sum_{n\leq N } q_n(x)
+\sum_{n>N} [p_n(x)-q_n(x)]\label{AE3}\\ &:=I_4+I_5+I_6,\nn
\end{align} where we choose $N$ to be the largest integer less than
$|x|^2/\log^2|x|$.

Note $$|I_4|\leq \P(\max_{n\leq N} |X_n|\geq |x|).$$ We estimate this using
truncation and Bernstein's inequality. Let $\xi_i =X_{i}-X_{i-1}$, define
$\xi'_i=\xi_i 1_{(|\xi_i|\leq N^{\frac12-\frac{\beta}4})}$, and
$X'_n=\sum_{i\leq n} \xi'_i$.  We have
\begin{align*}
\P(X_n\ne X'_n\mbox{ for some }n\leq N)&\leq
\P(\xi_n\ne \xi'_n\mbox{ for some }n\leq N)\\    &\leq N \max_{n\leq N}
\P(\xi_n\ne \xi'_n)\\ &\leq N c_1N^{(\frac12-\frac{\beta}4)(3+2\beta)}\leq
c_1N^{-\frac12-\frac{\beta}8}.
\end{align*} With our choice of $N$ we see that
\begin{equation}\label{AE4}
\P(X_n\ne X'_n\mbox{ for some }n\leq N)=o(1/|x|).
\end{equation} By Bernstein's inequality (\cite{Bennett})
\begin{align}
\P(\max_{n\leq N} |X'_n|\geq |x|)&\leq 2\exp\Big(-\frac{|x|^2} {2c_2N+
\frac23 |x| N^{\frac12 -\frac{\beta}{4}}}\Big)
\label{AE5}\\ &\leq 2e^{-c_3 \log^2|x|}=o(1/|x|).\nn
\end{align} Combining (\ref{AE4}) and (\ref{AE5}) yields the required bound
on $|I_4|$.

We can show $I_5=o(1/|x|)$ by straightforward estimates. Finally,  by
Proposition \ref{AP1},
$$|I_6|\leq \sum_{n>N} c_4 n^{-\frac32-\beta}=O(N^{-\frac12-\beta})
=o(1/|x|).$$ Summing the estimates for $I_4, I_5$, and $I_6$ shows
$I_3=o(1/|x|)$ and completes the proof.
\qed

The following result holds for all mean zero finite variance random walks in
any dimension $d$.  To keep the notation uniform we use $D(0,n)$ to denote
the  ball (if $d\geq 3$) or disc (if $d=2$) of radius
$n$ centered at the origin. When
$d=1$ we let $D(0,n)=( -n,n)$.

\begin{lemma}\label{ALesc} For some $c<\ff$
\begin{equation}
\E^x(T_{D(0,n)^c})\leq c n^2, \qquad x\in D(0,n),\hspace{ .2in} n\geq
1.\label{ALesc.1}
\end{equation}
\end{lemma}

\proof Let $T=\min\{j: |X_j-X_0|>2n\}$. By the invariance principle
\begin{eqnarray}
\P^x(T>c_1n^2)&=&\P^x(\sup_{j\leq c_1n^2} |X_j-X_0|\leq
2n)\label{Lesc.3}\\ &=&\P^0(\sup_{j\leq c_1n^2} |X_j-X_0|\leq 2n)\nn\\
&\leq &\rho<1\nn
\end{eqnarray} for all $x$ if we take $c_1=1$  and $n$ is large enough. Taking
$c_1$ larger if necessary, we get the inequality for all
$n$. Then letting $\theta_j$ be the usual shift operators and using the strong
Markov property
\begin{eqnarray}
\P^x(T>c_1(k+1)n^2)&\leq &\P^x(T\circ \theta_{c_1kn^2}> c_1n^2,
T>c_1kn^2)\label{Lesc.4}\\ &=&\E^x\Big[\P^{X_{c_1kn^2}}(T>c_1n^2);
T>c_1kn^2\Big]\nn\\ &\leq &\rho \P^x(T>c_1 kn^2).\nn
\end{eqnarray}
       Using induction
\begin{equation}
\P^x(T>c_1 kn^2)\leq \rho^k,\label{Lesc.10}
\end{equation}
  and our result follows easily.
\qed

Equation (6) of \cite{L2} does the simple random walk case of the following.

\begin{lemma}\label{LRatio} We have
$$\lim_{n\to \infty} \frac{\log T_{D(0,n)^c}}{\log n}=2, \qquad
\mbox{\rm
$\P^0$-a.s.}$$
\end{lemma}

\noindent{\bf Proof of Lemma \ref{LRatio}:} Let $\eps>0$. By Chebyshev and
Lemma \ref{ALesc},
$$\P^0(T_{D(0,n)^c}> n^{2+\eps})\leq \frac{\E^0 T_{D(0,n)^c}}{n^{2+\eps}}
\leq cn^{-\eps}.$$ So by Borel-Cantelli there exists $M_0(\omega)$ such that if
$m\geq M_0$, then $T_{D(0,2^m)^c}\leq (2^m)^{2+\eps}$.  If
$m\geq M_0$ and
$2^m\leq n\leq 2^{m+1}$, then
$$T_{D(0,n)^c}\leq T_{D(0,2^{m+1})^c} \leq (2^{m+1})^{2+\eps}\leq
2^{2+\eps} n^{2+\eps},$$ which, since
$\eps$ is arbitrary, proves the upper bound.

By Kolmogorov's inequality  applied to each component of the random walk,
$$\P^0(T_{D(0,n)^c}< n^{2-\eps}) =\P^0(\sup_{k\leq n^{2-\eps}} |X_k|>n)
\leq c\frac{\E^0 |X_{n^{2-\eps}}|^2}{n^2} \leq cn^{-\eps}.$$ So by
Borel-Cantelli there exists $M_1(\omega)$ such that if
$m\geq M_1$, then $T_{D(0,2^m)^c}\geq (2^m)^{2-\eps}$. If
$m\geq M_1$ and $2^m\leq n\leq 2^{m+1}$, then
$$T_{D(0,n)^c}\geq T_{D(0,2^m)^c}\geq (2^m)^{2-\eps}
\geq 2^{\eps-2} n^{2-\eps},$$ which proves the lower bound.
\qed

\vspace{2mm}
\noindent{\bf Acknowledgments} We would like to thank Greg Lawler for
making available to us some unpublished notes; these were very useful in
preparing Sections 2 and 9.

\bigskip
\noindent
\begin{tabular}{lll}  & Richard Bass& Jay Rosen\\
            & Department of Mathematics &  Department of Mathematics
\\
           &University of Connecticut & College of Staten Island, CUNY\\
           &Storrs, CT 06269-3009 & Staten Island, NY 10314\\
&bass@math.uconn.edu&jrosen3@earthlink.net
\end{tabular}

\end{document}